\documentclass[12pt, uplatex]{amsart}
\usepackage{amsmath,amsfonts,amsthm,amscd,amssymb,mathrsfs,amssymb,bm,pb-diagram, here}
\usepackage{setspace}
\usepackage[all,cmtip]{xy}

\usepackage{blkarray}
\usepackage{multirow}
\usepackage{mathtools}
\usepackage{graphicx}

\usepackage{mathrsfs, arydshln}
\usepackage{graphicx}
\usepackage[all]{xy}
\newtheorem{theorem}{Theorem}

\newtheorem{proposition}[theorem]{Proposition}
\newtheorem{lemma}[theorem]{Lemma}
\newtheorem{definition}[theorem]{Definition}

\newtheorem{corollary}[theorem]{Corollary}

\newtheorem{remark}[theorem]{Remark}
\newcommand{\aaa}{\alpha}
\newcommand{\bbb}{\beta}
\newcommand{\ccc}{\gamma}

\newcommand{\ddd}{\delta}
\newcommand{\DDD}{\Delta}

\newcommand{\lmd}{\lambda}

\newcommand{\CP}{\mathbb{CP}}

\newcommand{\CC}{\mathbb{C}}

\newcommand{\rank}{{\rm{rank}}}
\newcommand{\Sing}{{\rm{Sing}\,}}

\newcommand{\qdr}{\CP_1\times\CP_1}

\newcommand{\ol}{\overline}

\newcommand{\lras}{\,\longrightarrow\,}

\newcommand{\set}{\,|\,}

\newcommand{\proofend}{\hfill$\square$}

\newcommand{\inv}{^{-1}}

\newcommand{\Bs}{{\rm{Bs}}}

\newcommand{\ms}{\mathscr}
\newcommand{\minus}{\backslash}
\newcommand{\ptl}{\partial}

\newcommand{\qandq}{\quad{\text{and}}\quad}

\newcommand{\reg}{{\rm{reg}}}
\newcommand{\cusp}{_{\rm{cusp}}}

\newcommand{\Image}{{\rm{Image}}}
\newcommand{\Hess}{{\rm{Hess}\,}}

\numberwithin{equation}{section}
\numberwithin{theorem}{section}

\begin{document}
\bibliographystyle{alpha} 
\title[]
{On the cuspidal locus in the dual varieties of Segre quartic surfaces}
\author{Nobuhiro Honda}
\address{Department of Mathematics, Tokyo Institute
of Technology, 2-12-1, O-okayama, Meguro, 152-8551, JAPAN}
\email{honda@math.titech.ac.jp}
\author{Ayato Minagawa}
\address{Internet Initiative Japan Inc. Iidabashi Grand Bloom 2-10-2 Fujimi, Chiyoda-ku, Tokyo 102-0071, JAPAN}
\email{benkyo-gatizei@akane.waseda.jp}

\thanks{The first author was partially supported by JSPS KAKENHI Grant 16H03932.
}
\begin{abstract}
Motivated by a kind of Penrose correspondence,
we investigate the space of hyperplane sections of Segre quartic surfaces which have an ordinary cusp.
We show that the space of such hyperplane sections is empty for two kinds of Segre surfaces, and it
is a connected surface for all other kinds of Segre surfaces.
We also show that when it is non-empty, the closure of the space is either birational to the surface itself or birational to a double covering
of the surface, whose branch divisor consists of some specific lines on the surface.
\end{abstract}

\maketitle

%\tableofcontents
%\setcounter{tocdepth}{1}

\section{Introduction}
The so-called Penrose correspondence gives
a relationship between conformal differential geometry 
and complex geometry,
and typically it in particular means that the parameter space of smooth rational curves
on a complex manifold (i.e.\,a twistor space) is equipped with a special geometric structure,
if the normal bundles of the rational curves are of a particular form.
%The complex manifold is called the twistor space.
The most famous one would be the case where the rational curve is in a complex 3-manifold with the normal bundle being isomorphic to $\ms O(1)^{\oplus 2}$, and under this situation the space of rational curves
is a 4-dimensional manifold, equipped with a self-dual conformal structure \cite{AHS, Pen76}.
Another interesting instance is established by Hitchin \cite{Hi82-1},
in which case the rational curve is in a complex surface with normal bundle $\ms O(2)$.
The space of rational curves is a 3-dimensional manifold equipped with an Einstein-Weyl structure, which may be regarded as a conformal version of Einstein metrics. The complex surface in which the rational curves move is called the {\em minitwistor space},
and the rational curves are called {\em minitwistor lines}.

In the article \cite{HN11}, we showed that the same result about the presence of Einstein-Weyl structure still holds even when 
the rational curves on a complex surface have ordinary nodes,
as long as the complete family of such nodal rational curves on the surface is 3-dimensional.
These nodal rational curves are still called minitwistor lines, and
the number of nodes is called the {\em genus} of 
a minitwistor space.
%This gives a lot of new examples of minitwistor spaces.
In \cite{H20}, we showed that, under a natural minimality condition, compact minitwistor spaces with 
genus one are exactly a special algebraic surfaces, classically called {\em Segre surfaces}.
In modern language, these surfaces are nothing but the anti-canonical models of weak del-Pezzo surfaces of degree four,
and all of them are complete intersections of two quadrics in $\CP_4$.
The nodal rational curves in a Segre surface are  
sections by hyperplanes which are tangent to the surface at one smooth point.
This means that the space of minitwistor lines,
namely the Einstein-Weyl space associated to the minitwistor space,
is a Zariski-open subset of the projectively dual variety of the Segre surface.
Here, an interesting point is that 
the Zariski-open subset can never be the entire dual variety,
because nodal curves can always be deformed in the surface,  into a reducible curve or a curve with an ordinary cusp.
This means that, the Einstein-Weyl space is always non-compact and admits a natural compactification as a projective variety.
%An interesting point here is that,
%nodal rational curves can always be deformed in the surface,  into a reducible curve or a rational curve with an ordinary cusp.
%In other words, the dual variety of a Segre surface has a particular
%analytic subset whose points do not correspond to rational curve
%with one node.
As in \cite{H20}, we call 2-dimensional components of the added locus
in the compactification as
{\em divisors at infinity} of the Einstein-Weyl space.
Note that degeneration into a cuspidal curve can happen 
only when the genus of a minitwistor space is positive.

In the article \cite{H20}, for any Segre surface, we studied divisors at infinity which parameterize reducible hyperplane sections of the surface.
In particular we determined all components of such divisors,
and showed that they are  isomorphic to either
a projective plane or a smooth quadratic surface.
Also, we showed that the dual variety of a Segre surface intersects itself
along these components in such a way that the dual variety has
ordinary nodes along generic points of these components.

A purpose of the present article is to investigate 
the space of hyperplane sections of any Segre surface,
which have an ordinary cusp at a smooth point of the surface.
The closure of this space constitutes a subvariety in the dual variety
of the Segre surface, and we call it the {\em cuspidal locus}
in the dual variety of the surface.
The details on the cuspidal locus will be discussed in 
Section \ref{ss:cuspidal} for general surfaces embedded
in projective spaces.
In Section \ref{ss:cuspidal2}, by using deformation theory
of curves with singularities, we show that 
{\em the cuspidal locus in the dual variety of any Segre surface is 
reduced and 2-dimensional} if it is non-empty (Theorem \ref{thm:cuspcusp}),
and that {\em the dual variety has ordinary cusps along
the cuspidal locus} (Theorem \ref{thm:cuspcusp2}).
Thus, the cuspidal locus always constitutes a divisor at
infinity of the Einstein-Weyl space.

In Section \ref{ss:nemt}, we prove some basic properties
on Segre surfaces that will be needed in the rest of this article.
In Section \ref{ss:D1}, we determine the structure of certain divisor
in the incidence variety of the Segre surface, which is closely 
related to the cuspidal locus.
In Section \ref{ss:hsf}, we identify hyperplanes which cut out cuspidal curves from Segre surfaces, by using the incidence variety.
Consequently, we find that the cuspidal locus
can have three kinds of structures, depending on 
the number of such pencils on the Segre surfaces.
Subsequently in Section \ref{ss:pdc}, we express the number of 
pencils of double conics in terms of Segre symbols for Segre surfaces.
As a result, it turns out that {\em the cuspidal locus for two kinds of Segre surfaces are empty} (Corollary \ref{cor:empty}).
Also we see that among sixteen kinds of Segre surfaces, 
{\em the cuspidal locus for seven kinds of Segre surfaces are naturally 
birational to the surfaces themselves} (Corollary \ref{cor:bir}).

For the remaining seven kinds of Segre surfaces, the cuspidal locus
turns out to be birational to a double covering over the surfaces whose branch
divisor consist of some lines on the surfaces.
In Section \ref{ss:line1}, we show that every line on a Segre surface
which does not pass any singularities of the surface
is a simple branch divisor of the double covering
(Proposition \ref{prop:line1}).
In Sections \ref{ss:line2} and \ref{ss:line3},
we obtain similar results for any other lines on Segre surfaces.

In Section \ref{ss:conc}, by using the results in Section \ref{ss:line1}, we show that {\em the cuspidal locus in the dual variety of any Segre surface
is an irreducible surface if it is not empty} (Corollary \ref{cor:irr}).
This is the main result in this article.
Next we determine which lines on Segre surfaces are really
branch divisors of the above structure of double covering (Proposition \ref{prop:br3}).
By using these results, we construct a surface which is birational 
to the cuspidal locus in the dual variety of a smooth Segre surface,
and find that the cuspidal locus is a surface of general type
(Proposition \ref{prop:gt}).
In Section \ref{ss:sing}, we investigate singularities of the cuspidal locus
in the dual varieties of Segre surfaces,
which arise from the above double covering structure, 
and show that {\em the cuspidal locus always has ordinary cusps along some conic in the dual variety (Proposition \ref{prop:DDD})}. 

In Appendix, we investigate structure of the cuspidal locus
for some Segre surfaces, by using explicit equations of the surfaces.
As a result, we can precisely identify some divisor 
in the incidence variety which is studied in Section \ref{ss:D1}.
Also, we identify all lines on Segre surfaces which are not branch divisor
of the above double covering. This completes identification of 
the branch divisor of the above double covering which was postponed
in Section \ref{ss:conc}.

In general, the singular locus of the dual variety of a projective variety 
is known to be highly singular, but they are far from being well-understood. See \cite[Section 10]{Tev} for information
about this topic.
Because Segre surfaces are varieties of low degree and are of 
small codimension, 
they seem to be a nice class of varieties to investigate the singular locus of the dual variety in detail.
Indeed, the fact that Segre surfaces are of degree four
is essential in our investigation in many places.
Combined with the results in \cite{H20}, 
the present results would give a fairly complete understanding 
of the singular locus of the dual varieties of Segre surfaces.

\medskip\noindent
{\bf Acknowledgement.} We would like to 
thank Claude LeBrun for suggesting that degeneration of nodal curves (i.e.\,minitwistor lines) into cuspidal curves in a minitwistor space with positive genus would be of interest, which led us to the present work.

\section{Cuspidal locus}\label{s:cl}
\subsection{Some generalities on the cuspidal locus 
in the dual variety}
\label{ss:cuspidal}
First we briefly recall basic facts on projectively dual varieties. 
For more details, see Tevelev's book \cite{Tev} for example.
Let $X\subset\CP_n$ be a non-degenerate irreducible projective variety,
and write $X_{\reg}$ for the locus of smooth points of $X$.
We say that a hyperplane $H\subset\CP_n$ is tangent to $X$
if $H$ includes the tangent space $T_pX$ at some point $p\in X_{\reg}$.
This condition implies that the hyperplane section $H|_X$ has
a singularity at the point $p$.
The dual variety of $X$, which is denoted by $X^*$, 
is defined as the closure of the set of hyperplanes in $\CP_n$ which are
tangent to $X$, where the closure is taken in the dual projective space $\CP_n^*$.
The dual variety $X^*$ is a subvariety in $\CP_n^*$ and 
it is a hypersurface if the variety $X$ is not ruled,
namely when there is some point of $X$ which is not passed by a line in $X$
\cite[Theorem 1.18]{Tev}.
If $X$ is not ruled, for a generic hyperplane $H\in X^*$,
the section $H|_X$ has exactly one ordinary node as its all singularity \cite[Section 2.1.1]{VoII}.

It is also useful to introduce the incidence variety
$I(X)\subset \CP_n\times\CP_n^*$. This is by definition
the closure of the locus in $\CP_n\times\CP_n^*$
formed by a pair $(p, H)$
of a point $p\in X_{\reg}$ and
a hyperplane $H$ containing $T_pX$.
This is a subvariety in $\CP_n\times\CP_n^*$.
Then the dual variety is nothing but the image of
$I(X)$ under the projection to $X^*$,
and there is a diagram 
\begin{align}\label{diagram:fac2}
 \xymatrix{ 
&I(X) \ar[dl]_{\pi_1} \ar[dr]^{\pi_2} &\\
X && X^*
 }
\end{align}
where  $\pi_1$ and $\pi_2$ are restrictions to $I(X)$ of 
the projections from $\CP_n\times\CP_n^*$ to the 
two factors respectively.
Over the smooth locus $X_{\reg}$, the projection $\pi_1$ is a projective space bundle.
Another projection $\pi_2$ is birational as long as
$\dim X^*=n-1$.

Next, in order to define the cuspidal locus
for a non-degenerate irreducible projective variety $X\subset\CP_n$,
we consider the following locally closed subset of $I(X)$:
\begin{multline*}
\big\{(p,H)\in I(X) \,\big|\, {\text{$p\in X_{\reg}$}}, {\text{and
the singularity of
$H|_X$ at $p$ is
not an ordinary node}}\big\}.
\end{multline*}
Taking the closure in $I(X)$ of this subset,
we obtain a subvariety in $I(X)$.
In \cite[\S 10.2.4]{Tev}, 
the image of this subvariety under the projection $\pi_2$ to $X^*$
is denoted by $X^*\cusp$ when $X$ is smooth
(in that case one does not need to take the closure),
and for the moment we use this notation.
From the definition, for any hyperplane $H\in X^*\cusp$,
the section $H|_X$ has a non-nodal singularity at 
a smooth point of $X$, or perhaps a singularity at 
a singular point of $X$.
The variety $X^*\cusp$ is not necessarily irreducible.
Moreover, even when $\dim X = 2$, a generic point of an irreducible component
of $X^*\cusp$ can correspond to a hyperplane section
whose singularity is not an ordinary cusp.
So in this article, when $\dim X = 2$, 
we call the {\em cuspidal locus} in the dual variety $X^*$
for the union of all components
of the subvariety $X^*\cusp$ whose generic point 
corresponds to a hyperplane section
which has an ordinary cusp at a smooth point
of $X$.
So under this definition, the cuspidal locus really parameterizes
hyperplane sections which have an ordinary cusp.

Next, we give a description of the cuspidal locus
in terms of local coordinates when $X$ is a surface in $\CP_4$.
It will be used to investigate 
structure of the cuspidal locus of Segre surfaces.
We remark that it is not difficult to generalize the following description to the case where no constraint is supposed on 
 dimension and codimension of the subvariety $X$.

Let $S\subset\CP_4$ be a non-degenerate irreducible projective surface,
and pick any $p\in S_{\reg}$.
%and let $T_pS$ the tangent space at $p$.
By choosing  suitable non-homogeneous coordinates $(x,y,z,w)$, 
%, such that $T_pS=\{z_1 = \dots = z_r =0\}$.
we may suppose that, in a neighborhood of the point $p$, 
the pair 
$(x,y)$ works as holomorphic
coordinates on $S$ in a neighborhood of $p$,
and that 
there exist two holomorphic functions
$$ 
F=F(x,y)\qandq G=G(x,y),
%F_1=F_1(z_{r+1},\dots, z_{n}),\dots, F_r=F_r(z_{r+1},\dots, z_{n})
$$
defined in a neighborhood of $p$,
such that, around $p$, the surface $S$ is defined by the equations
$$
z=F(x,y),\quad w=G(x,y).
%z_1=F_1(z_{r+1},\dots, z_{n}),\dots, z_r=F_r(z_{r+1},\dots, z_{n}).
$$
Under these choices, for a point $q=(a,b)\in S_{\reg}$
in the neighborhood,
any hyperplane which contains the tangent plane $T_qS$ is of the form 
\begin{multline}\label{tpx2}
\lmd\left\{z - F(q)
- 
F_x(q) (x - a)
- %\frac{\ptl F}{\ptl y}
F_y(q) (y - b)\right\}\\
+
\mu
\left\{
w - G(q)
- 
G_x(q) (x - a)
- %\frac{\ptl G}{\ptl y}
G_y(q) (y - b)
\right\}=0
\end{multline}
for some $(\lmd,\mu)\in\CC^2\minus\{(0,0)\}$.
The pair $(\lmd,\mu)$ can be regarded as homogeneous
coordinates on the fiber of the projection $\pi_1:I(S)\to S$,
which is a projective line $\CP_{1}$.
Substituting the equations $z = F(x,y)$ and $w = G(x,y)$ into \eqref{tpx2},
we obtain a defining function of 
the intersection of the hyperplane \eqref{tpx2} with $S$.
This intersection has a non-nodal singularity at the point $q$
iff the Hessian (i.e.\,the determinant of the Hessian matrix) of the defining function vanishes at $q$.
It is immediate to see from \eqref{tpx2} that this Hessian is equal to the Hessian
$${\bf H} := \Hess \left(
\lmd F + \mu G
%\sum_{i=1}^r \lmd_i F_i(z_{r+1},\dots, z_n)
\right).
$$
Thus, to each $(\lmd,\mu)\in\CC^2\minus\{(0,0)\}$,
the section of $S$ by 
the hyperplane \eqref{tpx2}
has a non-nodal singularity at $q\in S_{\reg}$
if and only if ${\bf H}(q) = 0$.

Since the Hessian matrix is of size $2\times 2$,
as a homogeneous equation in $(\lmd,\mu)$,
the function $\bf H$ is quadratic in $\lmd$ and $\mu$.
More explicitly, we easily have
\begin{align}\label{Hess15}
{\bf H}= \Hess(F)\lmd^2
+ \big(F_{xx}G_{yy} + G_{xx}F_{yy} - 2F_{xy}G_{xy}\big)\lmd\mu
+ 
\Hess(G)\mu^2.
\end{align}
%As an equation of $(x,y)$,
The condition ${\bf H}=0$ is independent
of a choice of the coordinates $(x,y,z,w)$.
Therefore, the equation ${\bf H}=0$ defines a locally closed
subset in the incidence variety $I(S)$.
Throughout this paper, we denote the closure
of it by the bold letter $\bm D$.
This is a divisor in $I(S)$.
If the equations
\begin{align}\label{3cef}
\Hess(F) = \Hess(G) = F_{xx}G_{yy} + G_{xx}F_{yy} - 2F_{xy}G_{xy}=0
\end{align}
holds at some point $q\in S_{\reg}$, then whole the fiber 
$\pi_1\inv(q)=\CP_1$ 
is included in the divisor $\bm D$.
So, if there exists a curve on $S$ along which \eqref{3cef}
holds, then the fibers over the curve 
constitute a component of the divisor $\bm D$.
This component can have a multiplicity greater than one
in general.
Throughout this article, we denote $\bm D_1$ for the sum of all these components.
These components are mapped to curves by the projection $\pi_1:I(S)\to S$.
Of course, $\bm D_1$ might be the zero divisor.
We write $\bm D_S:=\bm D-\bm D_1$, so that 
$$
\bm D = \bm D_1 + \bm D_S
$$
holds. We always have $\bm D_S\neq 0$, and since the equation
\eqref{Hess15} is quadratic in $\lmd$ and $\mu$, 
the restriction of the projection $\pi_1:I(S)\to S$
to $\bm D_S$ is of degree two.
But this restriction is still not necessarily a finite morphism, because
there might exist an isolated point $q\in S_{\reg}$ such that 
the degeneracy condition \eqref{3cef} holds,
and also because $\bm D_S$ might include a curve in a fiber over some singular point of $S$.
Thus, what we can say in general is that the projection 
$\pi_1|_{\bm D_S}:\bm D_S\to S$ is a generically finite double covering over $S$.
The divisor $\bm D_S$ can be reducible, and in that case,
each of the two components is birational to $S$ by the projection $\pi_1$.
Regardless of whether $\bm D_S$ is irreducible or not,
we think $\bm D_S$ as a double covering over $S$.
If $\bm D_1 = 0$, then
the branch divisor of $\bm D_S\to S$ is the zero divisor of the discriminant of the quadratic polynomial \eqref{Hess15}.
If $\bm D_1 \neq 0$,
in order to obtain the equation of the branch divisor
of $\bm D_S\to S$,
we have to divided the quadratic polynomial \eqref{Hess15} by
a defining equation of $\bm D_1$,
before taking the discriminant.

Mapping the divisor $\bm D = \bm D_1+ \bm D_S$ in $I(S)$ to
the dual variety $S^*$ by the projection $\pi_2:I(S)\to S^*$,
we obtain a subvariety in $S^*$.
This is nothing but $S^*\cusp$ in the notation of \cite{Tev}
we have mentioned above.
This subvariety can be reducible in general, and 
the cuspidal locus in $S^*$ under current definition is some component of this subvariety
when it is non-empty.
A component of the image $\pi_2(\bm D)$
is not necessarily a component of the cuspidal locus in $S^*$ because a generic point
of the component can correspond to a non-cuspidal hyperplane section.

We end this subsection by showing that a straight-line on
a surface $S$ in $\CP_4$ is always a candidate of a branch divisor of
the generically finite double covering $\bm D_S\to S$:

\begin{proposition}\label{prop:line123}
If a non-degenerate irreducible surface $S\subset\CP_4$ contains a straight-line $l$
which is not contained in the singular locus of $S$, 
then the discriminant of the quadratic polynomial \eqref{Hess15} 
vanishes along $l$.
\end{proposition}

\proof
The discriminant of the polynomial $\bf H$ in \eqref{Hess15} can be 
easily calculated, and it can be written as
\begin{align}\label{disc}
 \big(F_{xx} G_{yy} - G_{xx} F_{yy} \big)^2
 + 4 \big(F_{xx}G_{xy} - F_{xy}G_{xx}\big)
 \big(F_{yy}G_{xy} - F_{xy}G_{yy}\big).
\end{align}
If $l$ is a line as in the proposition and 
$p\in l\cap S_{\reg}$, then
we can always take the non-homogeneous linear coordinates $(x,y,z,w)$
in such a way that $T_pS= \{z=w=0\}$
and $l=\{y=z=w=0\}$ hold.
These imply $F|_l = G|_l=0$,
which mean 
$$
\frac{\ptl^k F}{\ptl x^k}(x,0,0,0) = 
\frac{\ptl^k G}{\ptl x^k}(x,0,0,0) = 0
\quad{\text{for any $k\ge0$}}.
$$
From these, the discriminant \eqref{disc}
vanishes along the line $l$.
\proofend

\medskip
From the above relationship between the branch divisor of 
the generically finite double covering $\bm D_S\to S$ and 
the discriminant of the quadratic polynomial \eqref{Hess15},
a line $l\subset S$ is a branch divisor of $\bm D_S\to S$,
provided that 
some of the three coefficient functions of \eqref{Hess15} do not vanish
identically on the line $l$.
If all the three functions vanish along the line $l$,
and $m\,(>0)$ is  the minimal vanishing order of these functions
along $l$, then $\pi_1\inv(l)$ is included in the divisor $\bm D_1$
by multiplicity precisely $m$.
Then the line $l$ is really a branch divisor of 
$\bm D_S\to S$
if and only if the discriminant of the quadratic polynomial
${\bf H}/y^m$ vanishes along $l$ in the coordinates
of the previous proof.
In Sections \ref{s:scl} and \ref{s:scl2}, we will use these to determine
whether lines on Segre surfaces are
really branch divisors of the covering $\bm D_S\to S$.

\subsection{Cuspidal locus in the dual variety of a Segre surface}
\label{ss:cuspidal2}
A non-degenerate irreducible 2-dimensional subvariety $S\subset\CP_4$
is called a {\em Segre quartic surface} or simply 
a {\em Segre surface} if it is of degree four in $\CP_4$,
and is not a cone over a quartic curve in $\CP_3$ nor
a projection of a quartic surface in $\CP_5$.
All of them are complete intersections of two quadrics.
The normal forms of a pair of quadratic equations on $\CP_4$ are known,
and by using them, Segre surfaces
can be classified into 16 types.
Any Segre surface has at most a finite number of singularities,
and all of them are rational double points.
Any Segre surface has a finite number of straight-lines on it.
For more details on Segre surfaces, see \cite[Chapter\,XIII,
Section\,10]{HP} and \cite[Section 8.6]{Dol}.
For a brief account, see \cite{H20}.

In this subsection, we prove two basic properties on
the cuspidal locus in the dual varieties of Segre surfaces.
The first one concerns smoothness and dimension:

\begin{theorem}\label{thm:cuspcusp}
Let $S\subset\CP_4$ be any Segre surface,
and $H$ a hyperplane not passing any singularity
of $S$.
If $H$ belongs to the cuspidal locus in $S^*$ and $H|_S$ 
indeed has a single ordinary cusp as its all singularity,
then $H|_S$ is a rational curve and the cuspidal locus in $S^*$ is smooth and 2-dimensional at the point $H\in S^*$.
\end{theorem}

Note that the existence of a hyperplane section as in the theorem
is not obvious at all. 
Indeed, we will see in Section \ref{ss:pdc} that some Segre
surfaces do not have such a hyperplane $H$.

The second property is about singularity of the dual variety
along the cuspidal locus.

\begin{theorem}\label{thm:cuspcusp2}
The dual variety $S^*$ of any Segre surface $S$
has ordinary cusps along smooth locus in the cuspidal locus
formed by 
hyperplanes appearing in Theorem \ref{thm:cuspcusp}.
\end{theorem}

\noindent
This is similar to the property that 
the dual variety $S^*$ of a Segre surface $S$ has ordinary nodes
along the locus whose points correspond to hyperplane sections
of $S$ which have two ordinary nodes 
\cite[Proposition 3.15]{H20}.

\medskip
\noindent{\em Proof of Theorem \ref{thm:cuspcusp}.}
%First we show the assertion about $S^*\cusp$ itself.
Let $H$ be as in the theorem and put $C=H|_S$. 
By Lefshetz theorem, the curve $C$ is connected.
Let $p\in C$ be the ordinary cusp of $C$. By assumption, the curve $C$ is smooth except $p$, 
and hence $C$ is irreducible.
Since the arithmetic genus of $C$ is easily seen to be one,
the presence of a single cusp means that $C$ is a rational curve.

We denote $\Omega_S$ and $\Omega_C$ for the sheaves of K\"ahler differentials on $S$ and $C$ respectively,
and $\ms I_{C}$ for the ideal sheaf of $C$ in $\ms O_S$.
Associated to the embedding $C\subset S$, there is a standard exact sequence
\begin{align}\label{N'1}
0 \lras \ms I_{C}/\ms I_C^2 \lras \Omega_S|_C\lras \Omega_C\lras 0.
\end{align}
By applying $\ms H\!om_{\ms O_C}(\,\cdot\,,\ms O_C)$ to this sequence, 
noting that $\ms Ext^1(\Omega_S|_C,\ms O_C)=0$ since 
the surface $S$ is supposed to be smooth at points on $C$ and
it implies that $\Omega_S|_C$ is locally free,
we obtain an exact sequence
$$
0 \lras \Theta_C \lras \Theta_S|_C \lras N_{C} \lras T^1_C \lras 0,
$$
where $\Theta$ is the tangent sheaf
$\ms H\!om_{\ms O}(\Omega,\ms O)$,
$N_C$ means the normal sheaf $[C]|_C$
of $C$ in $S$,
and $T^1_C = \ms Ext^1(\Omega_C,\ms O_C)$.
Then the equi-singular normal sheaf $N'_C$ of $C$ in $S$ is defined 
as the kernel sheaf of 
the surjective homomorphism $N_{C} \lras T^1_C$
in this exact sequence.
In particular, we have a short exact sequence
\begin{align}\label{eqnn}
%0\lras T^0_C\lras T_S^0|_C \stackrel{\aaa}{\lras} N'_{C/S} \lras 0,\\
0\lras N'_{C} \lras N_{C} \lras T^1_C \lras 0.
\end{align}

In the following, for simplicity, we write $N$ and $N'$ for $N_C$ and $N'_C$ respectively.
By \cite{W74} (see also \cite[Proposition 1.1.9]{SeBook} and \cite[Section 1]{Tan}),  equi-singular displacements of $C$ in $S$
are governed by the cohomology groups of the sheaf $N'$.
More precisely, the Zariski tangent space at the point $C$  of the space of such displacements is identified with $H^0(N')$,
and an obstruction for smoothness of the last space at the point $C$ is in $H^1(N')$.
In particular, if $H^1(N')=0$, the space of such displacements is 
smooth at the point $C$, with the dimension being equal to $h^0(N')$.

Let $\ms J$ be the ideal sheaf on $S$ which is locally generated by the derivatives $\ptl f/\ptl x$ and $\ptl f/\ptl y$,
where $(x,y)$ are local coordinates on $S$ and $f=f(x,y)$
is a local defining function of the cuspidal curve $C$.
Note that we are assuming $C\cap\Sing S = \emptyset$.
Further, we put $\ms J_C:=\ms J\otimes_{\ms O_C}\ms O_C$.
These are called the Jacobian ideal sheaves
of the singular curve $C$.
As in \cite[p.\,111]{Tan} and \cite[\S 4.7.1]{SeBook}, 
the equi-singular normal sheaf satisfies
\begin{align}\label{NJ}
N' \simeq N\otimes_{\ms O_C} \ms J_C.
\end{align}
The way for calculating cohomology groups $H^q(N')$ using \eqref{NJ}
as well as the normalization of $C$ is briefly described in \cite[\S 4.7.1]{SeBook} for plane curves with ordinary nodes 
or cusps,
but we write it here with some detail since in Section \ref{ss:hsf}
we will need to calculate $H^q(N')$ for a curve which has
an $A_3$-singularity (tacnode).

To be explicit, let $\mu:S' \to S$ be the blowup of $S$ at
the cusp $p$ of $C$.
Then the restriction of $\mu$ to the strict transform $\tilde C$ of $C$ gives the normalization $\nu:\tilde C\lras C$.
As seen above, $\tilde C\simeq\CP_1$.
Take coordinates $(x,y)$ in a neighborhood of the cusp $p$
such that 
$C$ is locally defined by the equation $y^2 - x^3 = 0$.
Putting $y=ux$ for a new coordinate $u$ on the exceptional curve
of $\mu$,
the map $\mu$ is locally given by $(x,u)\mapsto (x,ux)$,
while $\tilde C$ is defined by $x=u^2$.
We also have $\ms J_p = (x^2, y)$ for the germ at $p$ of the Jacobian ideal sheaf, and therefore 
$(\mu^*\ms J)_{\tilde p} = (x^2,ux)$
for the germ at the point $\tilde p := \nu\inv(p)$.
Restricting $\mu^*\ms J$ to $\tilde C$ means a substitution
of $x=u^2$, so we obtain 
$$\big(\nu^*\ms J_C\big)_{\tilde p}
\simeq \Big(\big(\mu^*\ms J\big)\otimes_{\ms O_{S'}}\ms O_{\tilde C}\Big)_{\tilde p}
\simeq (u^4, u^3)
= (u^3).$$
This means 
\begin{align}\label{nuJ}
\nu^*\ms J_C \simeq \ms O_{\tilde C}(-3\tilde p).
\end{align}
Therefore, from \eqref{NJ} and $C^2 = 4$, we obtain
\begin{align}\label{NJ2}
\nu^*N' &\simeq \nu^*N\otimes_{\ms O_{\tilde C}} \nu^*\ms J_C\\
&\simeq \nu^*N \otimes_{\ms O_{\tilde C}} \ms O_{\tilde C}(-3\tilde p)\notag\\
&\simeq \ms O_{\tilde C}(1).\label{NJ3}
\end{align}
Taking the direct image of the isomorphism \eqref{NJ2}, 
from projection formula (which is available
since the sheaf $N$ is invertible),
we obtain $\nu_*\nu^*N'\simeq N\otimes \nu_*\nu^*\ms J_C$.
On the other hand, since the cusp has a single branch at $p$,
we have 
$\nu_*\nu^*\ms J_C \simeq \ms J_C$.
From these, we obtain $\nu_*\nu^*N'\simeq N\otimes \ms J_C
\simeq N'$.
By Leray spectral sequence, since we have
$R^q\nu_*(\nu^*N')=0$ for any $q>0$ as $\nu^*N'$ is invertible, 
it holds
$H^q(\tilde C, \nu^*N')\simeq H^q (C, \nu_*\nu^*N')$
for any $q\ge0$.
Thus we obtain 
$H^q(\tilde C, \nu^*N')\simeq H^q (C, N')$
for any $q\ge0$.
Hence from \eqref{NJ3}, we obtain 
$H^q(C, N') \simeq H^q(\tilde C, \ms O(1))$ for any $q\ge 0$.
This implies $H^1(C, N')=0$ and $h^0(C, N') = 2$.
These mean that equi-singular displacements of the cuspidal curve
$C$ in $S$ are parameterized by a smooth complex surface,
and this is what we need to show.
\proofend

\medskip
For the proof Theorem \ref{thm:cuspcusp2},
we recall basic results about deformations of $A_n$-singularity of a curve.
As is well-known, the versal family of $A_n$-singularity
is smooth, $n$-dimensional, and if we express the singularity by 
the equation $y^2 = x^{n+1}$ in $\CC^2$,  the versal family is concretely
given by 
\begin{align}\label{Kf1}
y^2 = x^{n+1} + s_1 x^{n-1} + s_2 x^{n-2} + 
\dots + s_n,
\quad s_1,\dots,s_n\in \CC.
\end{align}
If $C$ denotes the singular curve $y^2 = x^{n+1}$ for the moment, 
the tangent space of the parameter space $\CC^n$ of this family at the origin
$(s_1,\dots,s_n) = (0,\dots, 0)$ is naturally identified with
the space $H^0(T^1_C)$.
If we take a Galois cover of the parameter space
$\CC^n$ by introducing new parameters $t_1,\dots, t_{n+1}$
by putting
\begin{align}\label{Gc}
s_1 = \sum_{1\le i<j\le n+1} t_it_j,
\quad
s_2 = \sum_{1\le i<j<k\le n+1} t_it_jt_k,\quad
\dots,\quad
s_{n} = t_1 t_2 \cdots t_{n+1}
\end{align}
and imposing the condition $t_1 + \dots + t_{n+1} = 0$, 
then as a base change of the family \eqref{Kf1}, we obtain another $n$-dimensional family 
\begin{align}\label{Gc1}
y ^2 = (x-t_1) (x-t_2) \cdots (x-t_{n+1}),  \quad
(t_1,\dots, t_{n+1})\in \CC^{n+1},\quad
\sum_{1\le i\le n+1} t_i = 0.
\end{align}
If $k$ is an integer satisfying $0<k<n$,
the locus of points whose fibers have 
an $A_k$-singularity as their all singularity is concretely given by 
the conditions that exactly $k$ among $t_1,\dots, t_{n+1}$ are equal
and that the remaining parameters are not equal.
In particular, when $n=2$, the fiber over the point $(t_1,t_2,t_3)$ 
with $t_1 + t_2 + t_3 = 0$ of the family \eqref{Gc1} has
an $A_1$-singularity iff $(t_1,t_2,t_3)$ is of the form
$(t,t,-2t),\,(t,-2t,t)$ or $(-2t,t,t)$ for some $t\neq 0$.
Hence, from \eqref{Gc}, the fiber over the point $(s_1,s_2)$ 
of the versal family \eqref{Kf1} has an $A_1$-singularity iff 
$(s_1,s_2) = (-3t^2, -2t^3)$ for some $t\neq 0$.
We denote the locus formed by these points in $\CC^2$ by $\ms A_1$.
Then $\ms A_1$ is locally closed, smooth, 1-dimensional, and 
its closure in $\CC^2$ has an $A_2$-singularity at the origin.
With these preliminaries, we provide:

\medskip\noindent
{\em Proof of Theorem \ref{thm:cuspcusp2}.}
We keep notations in the proof of Theorem \ref{thm:cuspcusp}.
We first show $H^1(N)=0$ and $H^0(N)\simeq\CC^4$.
Let $C$ be any hyperplane section of the Segre surface $S$
which has a single ordinary cusp at a smooth point of $S$ as its all singularity.
By using that the surface $S\subset\CP_4$ is a complete intersection 
of two quadrics, we readily obtain that $C$ belongs to the 
anti-canonical class on $S$.
This means that there is
an exact sequence
\begin{align}\label{ses03}
0 \lras \ms O_S \lras K_S\inv \lras N\lras 0.
\end{align}
This time we mean by $\mu:\tilde S\to S$ the minimal resolution of 
all singularities of $S$.
Since all the singularities of $S$ are rational double points
(\cite[Theorem 8.1.11]{Dol}),
we have $K_{\tilde S}\simeq \mu^*K_S$.
Further, %since $\tilde S$ is a weak del-Pezzo surface,
$K_{\tilde S}\inv$ is nef and big (\cite[p.\,355]{Dol}).
Furthermore, by Leray spectral sequence, for any 
invertible sheaf $\ms L$ on  $S$, 
we have $H^q(\tilde S, \mu^*\ms L)\simeq
H^q(S,\ms L)$ for any $q\ge0$.
From these and Kodaira-Ramanujan vanishing theorem,
it readily follows that $H^q(K_S\inv)=H^q(\ms O_S) = 0$
for any $q>0$.
Hence from the cohomology exact sequence of \eqref{ses03},
we obtain $H^1(N) = 0$.
Also, as showed in the proof of \cite[Lemma 3.8]{H20},
the restriction homomorphism
\begin{align}\label{isom0}
H^0\big(\CP_4,\ms O(1)\big) \lras 
H^0\big(K_S\inv\big) 
\end{align}
is isomorphic, and hence $h^0(K_S\inv) = 5$.
So again from the cohomology exact sequence of \eqref{ses03},
we obtain $h^0(N) = 4$.

The fact $H^1(N) = 0$ means that any first order displacement of $C$ in the surface $S$ (without any constraint for displacements this time) is unobstructed,
and the tangent space at the point $C$ of the parameter space of the versal family of displacements of $C$ in $S$  is identified with $H^0(N)\simeq\CC^4$.

From the short exact sequence \eqref{eqnn}, 
since $H^1(N')=0$ as we have shown
in the proof of Theorem \ref{thm:cuspcusp}, we obtain an exact sequence
\begin{align}\label{eqnn2}
%0\lras T^0_C\lras T_S^0|_C \stackrel{\aaa}{\lras} N'_{C/S} \lras 0,\\
0\lras H^0\big(N'\big) \lras H^0\big(N\big) \lras 
H^0\big(T^1_C\big) \lras 0.
\end{align}
Let $H\subset\CP_4$ be the hyperplane which satisfies 
$H|_S = C$, and
$B$ a neighborhood of the point $H$ in the dual space $\CP^*_4$.
To each point of $B$, we can naturally associate a hyperplane section of $S$,
so $B$ can be regarded as a parameter space of 
displacements of $C$ in $S$. 
As above, we have a natural isomorphism $T_H B \simeq H^0(N)$.
By versality, after shrinking $B$ if necessary, 
there is an induced holomorphic map
$f:B\to \CC^2$ from $B$ to the 
parameter space of the versal family of an $A_2$-singularity,
while the differential $df:T_H B\to T_0\CC^2$
is identified with the map
 $H^0\big(N\big) \lras 
H^0\big(T^1_C\big)$ in \eqref{eqnn2}.
Writing $\ms A_1\subset\CC^2$ for the locus 
formed by points whose fibers have $A_1$-singularity as 
presented right before the present proof,
because a generic point of the dual variety $S^*$ 
corresponds to a hyperplane section which has an $A_1$-singularity,
we have 
$$S^*\cap B=f\inv\big(\ms A_1\cup\{H\}\big).$$
Since the map $H^0\big(N\big) \lras 
H^0\big(T^1_C\big)$ in \eqref{eqnn2} is surjective, 
%by the implicit function theorem, 
again shrinking the neighborhood $B$ if necessary, 
there are coordinates $(z_1, z_2, z _3, z_4)$ on $B$
such that $f$ takes the form 
$(z_1, z_2, z _3, z_4)\longmapsto (z_1, z_2)$.
As seen right before the present proof,
the locus $\ms A_1\cup\{H\}$ is 1-dimensional and 
has an ordinary cusp at the point $H$ as its only singularity.
These imply that the dual variety $S^*$ has ordinary cusps along the smooth surface $f\inv(0,0)=\{z_1=z_2=0\}$ in $B$,
formed by hyperplane sections which have ordinary cusp 
as its all singularity.
\proofend

\section{The structure of the cuspidal locus}
\subsection{Some general properties of Segre surfaces}
\label{ss:nemt}
In this subsection, we prove a few basic properties
of Segre surfaces that will be used throughout the rest of this article.
%Next, by using these properties,  
%we show that the cuspidal locus in the dual variety
%of any Segre surface is not empty.
We begin with an easy but useful property about curves
with low degree on Segre surfaces.

\begin{proposition}\label{prop:SI1}
Let $S\subset\CP_4$ be a Segre surface,
and $\mu:\tilde S \to S$ the minimal resolution of all singularities of $S$.
The self-intersection number in $\tilde S$ of the strict transform of a line, an irreducible conic, or a rational 
normal curve of degree three, lying on $S$, is 
$(-1)$, $0$ and $1$ respectively.
\end{proposition}

Of course, when the curve (a line, a conic, or a rational normal curve of degree three) does not pass any singularity of $S$,
this implies that the self-intersection number of the curve in $S$ is 
$(-1)$, 0 and $1$ respectively.
We note that the proposition in particular means that an irreducible conic on any Segre surface $S$ is always a member of a pencil of conics on $S$,
and the pencil can have a base point only at a singularity
of $S$.

\medskip\noindent
{\em Proof of Proposition \ref{prop:SI1}.}
We write $C$ for the curve in $S$ and $\tilde C$ for the 
strict transform of $C$ into $\tilde S$.
As remarked in the proof of Theorem \ref{thm:cuspcusp2}.
if $H$ is the hyperplane class on $\CP_4$,
we have $H|_S = K_S\inv$ and $K_{\tilde S}\simeq \mu^*K_S$.
Noting that the restriction $\mu|_{\tilde C}$ gives an isomorphism $\tilde C\simeq C$ as $C$ is smooth,
we have
\begin{align*}
K_{\tilde C} &\simeq K_{\tilde S}|_{\tilde C} \otimes N_{\tilde C/\tilde S}\\
&\simeq (\mu^*K_{S})|_{\tilde C}\otimes N_{\tilde C/\tilde S}\\
&\simeq \mu^*\big(K_{S}|_{C}\big)\otimes N_{\tilde C/\tilde S}\\
&\simeq (\mu|_{\tilde C})^*\big(-H|_{ C}\big)\otimes N_{\tilde C/\tilde S}.
\end{align*}
Hence, since $\tilde C$ is a smooth rational curve in the present situation, 
comparing the degrees of both sides, we obtain
$$
-2 = -H.\,C + \tilde C^2.
$$
From this we readily obtain the assertions of the proposition.
\proofend

\medskip
We also have 
\begin{proposition}\label{prop:SI2}
Any Segre surface $S$ does not have a reduced cubic curve which is contained in a 2-plane in $\CP_4$.
\end{proposition}

\proof
Let $C\subset S$ be such a cubic curve, and $P$ the 2-plane
containing $C$.
Then for any hyperplane $H$ containing $P$,
we have $H|_S = C + l$ for a line $l$.
Since we have an isomorphism \eqref{isom0},
the line $l$ moves as $H$ moves.
This contradicts finiteness of lines on $S$.
Hence, such a curve $C$ does not exist.
\proofend

\medskip
In the following, by a tangent space of a Segre surface at a smooth point, we always mean a closed linear subspace
in the projective space $\CP_4$. 
So it is a projective 2-plane, not
just a vector space.
In general, if $l$ is a line lying on a subvariety, then clearly the tangent space
at any smooth point belonging to $l$ contains $l$.
The next lemma says that for any Segre surface $S$,  
lines on $S$ are all curves on $S$
which are contained in a tangent plane at some smooth point of 
$S$.

\begin{lemma}\label{lemma:tan1}
Let  $S\subset\CP_4$ be a Segre surface,
$p$ a smooth point of $S$,
and $T_pS$ the tangent plane in the above sense.
If the intersection $S\cap T_pS$ contains a curve,
then any of its irreducible component is a line through the point $p$.
\end{lemma}

For the proof of this lemma, it is convenient to introduce the following
notations. They will be used throughout the rest of this article.
\begin{definition}{\em
For a smooth point $p$ of a Segre surface $S\subset\CP_4$,
we denote $(T_pS)^*$
%$\langle T_pS \rangle$ 
for the pencil of hyperplanes in $\CP_4$
which contain the tangent plane $T_pS$.
Also, we write  $(T_pS)^*|_S$
%$\langle T_pS \rangle|_S$ 
for the pencil on $S$
whose members are of the form $H|_S$ for some $H\in 
(T_pS)^*$. \proofend
}
\end{definition}

\noindent
{\em Proof of Lemma \ref{lemma:tan1}.} 
Since $\deg S = 4$, the intersection $S\cap T_pS$ can contain
a curve whose degree is at most four.
If $S\cap T_pS$ contains an irreducible quartic curve $C$,
then by the same reason, any hyperplane $H\in (T_pS)^*$ satisfies $S\cap H= C$. In particular, the intersection $S\cap H$ is 
independent of a choice of such a hyperplane $H$.
This contradicts the basic isomorphism \eqref{isom0}.
%\begin{align}\label{isom2}
%H^0\big(\CP_4,\ms O(1)\big) \lras 
%H^0\big(K_S\inv\big) 
%\end{align}
%is isomorphic. (See \eqref{isom0} in the last subsection.)
Hence, $S\cap T_pS$ does not contain an irreducible quartic curve.
Next, if $S\cap T_pS$ contains an irreducible cubic curve $C$,
for any hyperplane $H$ containing $T_pS$,
we have  $S\cap H= C+l$ for some line $l$.
%there exists a line $l$ such that  $S\cap H= C+l$ holds.
Again from the isomorphism \eqref{isom0},
if we move $H$ in the pencil $(T_pS)^*$,
the line $l$ also really moves.
This contradicts finiteness of the number of lines on any
Segre surfaces.
Hence $S\cap T_pS$ does not contain an irreducible cubic curve.

Next, suppose that $S\cap T_pS$ contains an irreducible conic $C$.
Then for any $H\in (T_pS)^*$,
$H|_S$ is of the form $C+D$ for some conic $D$ which is possibly
reducible or non-reduced.
Again from the isomorphism \eqref{isom0},
this conic $D$ has to move as $H$ moves in $(T_pS)^*$.
Hence, again from the finiteness of lines, the conic $D$ has to be irreducible for a generic $H\in (T_pS)^*$.
But then $D$ would also be smooth, and therefore
$D$ has to pass the point $p$ since the curve $C+D$ is singular at $p$ as $H\in (T_pS)^*$.
This means that 
there is a pencil of conics on $S$ which has a base point
belonging to $S_{\reg}$.
This contradicts Proposition \ref{prop:SI1}.
(See the comment right after the proposition.)
Thus, the intersection $S\cap T_pS$ does not contain an irreducible conic.

Therefore, any curve  in $S\cap T_pS$ is a line.
Let $l$ be such a line and suppose that $p\not\in l$.
If $H\in  (T_pS)^*$, then 
$H|_S = l + D$ for some cubic $D$ which is possibly
reducible or non-reduced.
Since $p\not\in l$, the curve $D$ has to have at least a double point at $p$.
We show that the pencil 
$$
(T_pS)^*|_S-l = \left\{H|_S - l\set 
H\in (T_pS)^*\right\}
$$ 
formed by the residual cubic curves has no fixed component.
Suppose not and let $F$ be the fixed component.
Then $F$ is a curve on the 2-plane $T_pS$,
so it has to be a line as above.
Moreover, $F$ goes through $p$, since otherwise
the movable part of the pencil $(T_pS)^*|_S-l-F$ consists of conics whose general members are irreducible,
and this contradicts that any member of 
the pencil $(T_pS)^*|_S$ is singular at $p$.
Thus, under the assumption on the presence of the fixed component $F$,
a generic member of $(T_pS)^*|_S$
can be written as $l + F + C$, where $F$ is a line through $p$
and $C$ is an irreducible conic through $p$.
This again contradicts Proposition \ref{prop:SI1},
and we obtain that the residual pencil $(T_pS)^*|_S-l$ has no fixed component.
From the finiteness of lines on $S$, this means that 
a generic member $D$ of this pencil is an irreducible cubic curve.
Because we are supposing $p\not\in l$, the cubic curve $D$ itself is singular at $p$.
Since $D$ is an irreducible cubic curve, this implies that $D$ lies on some 2-plane $P$
since otherwise $D$ would be a rational normal curve, which is smooth.
This contradicts Proposition \ref{prop:SI2}.
Therefore, the point $p$ belongs to $l$, as desired.
\proofend

\medskip
By using this lemma,  we next show that there is no point $p\in S_{\reg}$ such that any member of the pencil $(T_pS)^*|_S$
has a triple point at $p$.

\begin{lemma}\label{lemma:dbl}
Let $S\subset\CP_4$ be a Segre surface and $p$ a smooth point of $S$.
Then the multiplicity of a generic member
of the pencil $(T_pS)^*|_S$ is precisely two.
\end{lemma}

\proof
First we show that if all members
of the pencil $(T_pS)^*|_S$ have a triple point at $p$,
then there has to exist a line on $S$ through $p$.
Let $\mu:S'\to S$ be a blowup at $p$,
and $E$ the exceptional curve.
If the hyperplane section $H|_S$ would have 
a triple point at $p$ for any  $H\in (T_pS)^*$, then
the linear system $|\mu^*H - 3E|$ on $S'$ is also a pencil.
Then since 
$
(\mu^*H - 3E)^2 = 4 - 9 = -5<0,
$
the pencil $|\mu^*H - 3E|$ would have a fixed component.
If $E$ would be a fixed component, then the system
$|\mu^*H - 4E|$ is a pencil,
and again from the self-intersection number, this has
to have a fixed component.
This process cannot continue infinitely many times.
So there exists some $m>2$ such that 
the system $|\mu^*H - mE|$ is a pencil and has a fixed component
other than $E$.
If we take the image of this component to $S$ by $\mu$, 
we obtain a curve on $S$, and it has to be contained in 
the 2-plane $T_pS$ since it has to be a base curve of the pencil 
$(T_pS)^*|_S$.
By Lemma \ref{lemma:tan1}, the last curve is a line through $p$.

Hence, if any member of the pencil $(T_pS)^*|_S$ has a triple point at $p$,
a generic member of the pencil is of the form $l + D$ for some line $l$ through $p$ and a cubic curve $D$ which may 
be reducible or non-reduced.
The curve $D$ has to be singular at $p$ since $l+D$ is assumed to 
have a triple point at $p$.
If a generic member of the residual pencil
$(T_pS)^*|_S - l$ of cubic curves 
would be reducible, by finiteness of lines on $S$,
the pencil has to have another line $l'$ as a fixed component.
If $p\not\in l'$,
then any member of the pencil $(T_pS)^*|_S - l - l'$ of conics has to have a double point at $p$, and this cannot happen by finiteness of lines on $S$.
Hence $p\in l'$. 
But still the pencil $(T_pS)^*|_S - l - l'$ of conics has 
$p$ as a base point,
which again contradicts Proposition \ref{prop:SI1}.
Hence a generic member of the pencil $(T_pS)^*|_S - l$ is irreducible.
So there exists an irreducible cubic curve $D$ on $S$
which is singular at $p$.
The curve $D$ has to lie on a 2-plane, and this contradicts
Proposition \ref{prop:SI2}.
Therefore, a generic member
of the pencil $(T_pS)^*|_S$ 
does not have a triple point at $p$.
\proofend

\subsection{Structure of the divisor $\bm D_1$}
\label{ss:D1}
As before, let $S\subset\CP_4$ be a Segre surface.
In this subsection, by using the results in the previous subsection,
we describe the structure of the divisor $\bm D_1$ 
in the incidence variety $I(S)$.
(See Section \ref{ss:cuspidal} for the definition of the divisor $\bm D_1$.)

Let $p$ be a smooth point of $S$,
$\mu:S'\to S$ a blowup at $p$, and $E$ the exceptional curve.
Then there is the following natural isomorphism between
pencils:
\begin{align}\label{TpSpen}
(T_pS)^*|_S\simeq |\mu^* H - 2E|.
\end{align}
If the point $p$ is on a line $l\subset S$, then
the strict transform $l'$ of $l$ into $S'$ is obviously 
a fixed component of the pencil $|\mu^* H - 2E|$.
The following proposition can be thought of as a kind of its converse.

\begin{lemma}\label{lemma:irr}
Let $S\subset\CP_4$ be a Segre surface
and $p$ a smooth point of $S$.
If the pencil $|\mu^*H - 2 E |$ on the blowup has a base point,
% belonging to
%the exceptional curve $E$, 
then there exists a line lying on $S$
which passes the point $p$.
\end{lemma}

\proof
It is enough to show that,
the pencil $|\mu^*H - 2 E |$ has a fixed component which is different from $E$,
since the image of such a component to $S$ by $\mu$ is a fixed component
of the pencil $(T_pS)^*|_S$, and therefore it is contained 
in the tangent plane $T_pS$, which has to be a line through $p$
by Lemma \ref{lemma:tan1}.

First, by Lemma \ref{lemma:dbl}, the exceptional curve $E$
of $\mu$ cannot be a base curve of the pencil $|\mu^*H - 2 E |$
because if so, any member of the pencil $(T_pS)^*|_S$
would have a triple point at $p$.
So assume that the pencil $|\mu^*H - 2 E |$ has no fixed component
and derive a contradiction. 
Since we are supposing $\Bs\,|\mu^*H - 2 E |\neq\emptyset$,
this assumption means that the pencil $|\mu^*H - 2 E |$ has an 
isolated fixed point.
Suppose that some (and hence a generic) member of
the pencil $|\mu^*H - 2 E |$ is irreducible.
Since the arithmetic genus of the line bundle $\mu^*H - 2 E$
is readily seen to be zero, this implies that 
the irreducible member of $|\mu^*H - 2 E |$ is smooth and rational.
We also obtain that the self-intersection number of the rational curve
is zero.
These mean that $\Bs\,|\mu^*H - 2 E |=\emptyset$.
Hence any member of the pencil $|\mu^*H - 2 E |$ is reducible.
Therefore, any member of the pencil $(T_pS)^*|_S$ is also reducible.
Hence, from the finiteness of lines on $S$,
the only possibility for the form of a generic member of the pencil $(T_pS)^*|_S$
is $C_1 + C_2$, where $C_1$ and $C_2$ are irreducible conics.
Since $C_1$ and $C_2$ are then smooth,
both of them must pass the point $p$.
Therefore at least one of the conics $C_1$ and $C_2$
moves on $S$ while passing the smooth point $p$ of $S$.
But this again contradicts Proposition \ref{prop:SI1}.
\proofend

\medskip
From this, it is not difficult to obtain the following

\begin{lemma}\label{lemma:bpf2}
Let $S\subset\CP_4$ be a Segre surface and $p$ a smooth point of $S$.
If there exists no line on $S$ through $p$, 
then the pencil $|\mu^* H - 2E|$ has exactly two members
whose restrictions to $E$ are double point.
\end{lemma}

\proof
It is elementary to prove that,
under a suitable choice of homogeneous coordinates  $(u,v)$ on $\CP_1$,
any 1-dimensional subsystem of the linear system
$|\ms O(2)|$ on $\CP_1$ is generated by
either $u^2$ and $v^2$, or $u^2$ and $uv$.
These are distinguished by presence of a base point of the subsystem.
By Lemma \ref{lemma:irr}, under the present assumption on the point $p$,
the pencil $|\mu^* H - 2E|$ is base point free.
So it does not have a base point also on $E$.
This means that the restriction homomorphism
\begin{align}\label{rE}
r_E:H^0(\mu^*H - 2E)\lras H^0(\ms O_E(2))
\end{align}
has a 2-dimensional image,
and that the image is generated by 
$u^2$ and $v^2$ for suitable homogeneous coordinates $(u,v)$ on $E$.
Obviously, the two members of the pencil 
$(T_pS)^*|_S$ which correspond to these generators
are all members
which satisfy the property in the lemma.
\proofend

\medskip
From this lemma, we obtain the following result about
the defining function $\bf H$ of the divisor $\bm D$ in $I(S)$.
(See Section \ref{ss:cuspidal} for the definition of 
$\bf H$ and $\bm D$.)

\begin{proposition}\label{prop:H}
Let  $S\subset\CP_4$ be a Segre surface and $p\in S_{\reg}$
a  point  which does not belong to any line on $S$.
Then we have ${\bf H}(p)\neq 0$ in the sense that 
the quadratic polynomial ${\bf H}(p)$ in $\lmd,\mu$ is not
identically zero.
Further, the quadratic equation ${\bf H}(p)=0$ does not have
a double root.
\end{proposition}

\proof
As in the previous proof,
under the assumption on the point $p$,
the restriction map 
$r_E$ as in \eqref{rE} has a two-dimensional image.
So not all elements of $\Image(r_E)$ can have a double root on $E$.
%and it is generated by $\lmd^2$ and $\mu^2$.
Because zeroes of ${\bf H}(p)$ correspond to hyperplanes
whose sections have non-nodal singularity at $p$,
this means that a generic member of the pencil 
$(T_pS)^*|_S$ has an ordinary node at $p$.
This implies ${\bf H}(p)\neq 0$.
Since each solution of the equation ${\bf H}(p) = 0$ about 
$(\lmd,\mu)\in\CP_1$ corresponds to a hyperplane
$H\in (T_pS)^*$ such that the singularity 
of $H|_S$ at $p$ is non-nodal, the equation ${\bf H}(p) = 0$ has a double root
iff such a hyperplane is unique.
On the other hand, via the blowup $\mu$ at $p$,
such a hyperplane gives
an element of the pencil $\big|\Image(r_E)\big|$ on $E$ which is a double point.
By Lemma \ref{lemma:bpf2},
when $p$ is not on a line on $S$, the last pencil has
exactly two such members. 
These mean that the equation ${\bf H}(p) = 0$ does not have a double root.
\proofend

\medskip As an immediate consequence, we can identify 
any component of the divisor $\bm D_1$ in $I(S)$ as follows.

\begin{corollary}\label{cor:D1}
For any Segre surface $S\subset \CP_4$,
any irreducible component of the divisor $\bm D_1$
is of the form $\pi_1\inv(l)$, where $l$ is a line
lying on $S$.
\end{corollary}

Note that this does not assert that $\pi_1\inv(l)$ 
is a component of the divisor $\bm D_1$ for any line $l$ on $S$.
Indeed, we will see that this is not correct in general.

From the corollary, we obtain the following information about
the cuspidal locus in the dual varieties of Segre surfaces.
\begin{proposition}\label{prop:cuspD2}
Let $S\subset\CP_4$ be any Segre surface.
Then the cuspidal locus in the dual variety $S^*$ is
contained in the image $\pi_2(\bm D_S)$, where $\pi_2:I(S)\to S^*$ is the projection introduced in Section \ref{ss:cuspidal}.
\end{proposition}

\proof
From our definition of the cuspidal locus given in Section \ref{ss:cuspidal},
any component of the cuspidal locus has a point $H$
such that $H|_S$ has an ordinary cusp at some point $p\in S_{\reg}$.
Then the point $(p,H)\in I(S)$ belongs to either $\bm D_1$
or $\bm D_S$.
If $(p,H)\in \bm D_1$, then from Corollary \ref{cor:D1},
the smooth point $p$ belongs to a line $l$ on $S$,
and the hyperplane $H$ includes $l$.
This means that the section $H|_S$ contains $l$ as a component.
Therefore, $H|_S$ does not have an ordinary cusp at $p$,
which contradicts our choice of $H$.
This implies $(p,H)\in \bm D_S$.
So the hyperplane $H$ belongs to the image $\pi_2(\bm D_S)$.
\proofend

\subsection{Hyperplane sections from the divisor $\bm D_S$}
\label{ss:hsf}
From Proposition \ref{prop:cuspD2}, for any Segre surface $S$, 
the investigation of the cuspidal locus in the dual variety
$S^*$ is reduced to that of the divisor
$\bm D_S$ in $I(S)$ and its image to $S^*$ by the projection 
$\pi_2:I(S)\to S^*$.
In this subsection, we first investigate the structure of hyperplane sections of a Segre surface which are obtained from generic points of the divisor $\bm D_S$.
Next, by using it, we show that there are three possibilities for the concrete forms of hyperplane sections which are obtained from generic points of $\bm D_S$, and then derive some conclusions about irreducibility
of the cuspidal locus in the dual varieties of Segre surfaces.

Take any point $p\in S_{\reg}$ which does not belong to any line on $S$.
By Proposition \ref{prop:H}, 
there exist precisely two hyperplanes belonging to 
the pencil $(T_pS)^*$
whose intersections with $S$ have non-nodal singularity at $p$.
We denote $H_1$ and $H_2$ for these hyperplanes.
We then have $(p,H_1)\in \bm D_S$ and $(p,H_2)\in \bm D_S$.

\begin{proposition}\label{prop:2mem}
Let $C$ be any one of the two hyperplane sections $H_1|_S$
and $H_2|_S$.
If $C$ is irreducible, then $C$ has an ordinary cusp at $p$
as its all singularity.
If $C$ is reducible, 
then exactly one of the following holds:
\begin{itemize}
\item%[\em (i)]
$C$ is an irreducible conic with multiplicity two,
\item%[\em (ii)]
$C$ consists of two irreducible conics
which are not co-planer and which are tangent to each other at $p$ with order two.
\end{itemize}
\end{proposition}

\proof
Since $p_a(C) = 1$ for the arithmetic genus of $C$,
when $C$ is irreducible,  it is a rational curve which has an
ordinary node or an ordinary cusp at $p$ as its all singularity.
But since the singularity is non-nodal as above, 
it has to be an ordinary cusp.

Next assume that $C$ is reducible and 
has some line $l$ as a component. 
Write $C = l + D$, where $D$ is a cubic curve.
Since $l$ cannot go through $p$ from the choice of $p$, 
the cubic $D$ has to have a non-nodal singularity at $p$.
From absence of lines through $p$,
this means that the cubic $D$ is irreducible and $p$
is an ordinary cusp of $D$.
Therefore, $D$ is contained in a 2-plane. 
This contradicts Proposition \ref{prop:SI2}.
Hence, the curve $C$ does not contain any line.
So we have $C= C_1+C_2$, where both $C_1$ and $C_2$ are
irreducible conics.
Since both are smooth, they pass the point $p$.
If $C_1 = C_2$, then $C$ is a double irreducible conic
as in the first item in the proposition.

Suppose $C_1\neq C_2$. They do not intersect
transversely at the point $p$ as $C$ is non-nodal at $p$.
Because of the isomorphism \eqref{isom0},
they cannot lie on the same 2-plane.
Since $C_1$ and $C_2$ are tangent at $p$,
this means $C_1\cap C_2 = \{p\}$.
Because $p$ is a smooth point of $S$ and 
$C_1+C_2$ is a hyperplane section of $S$,
this implies that 
the smooth curves $C_1$ and $C_2$ do not pass any singularity of $S$. From Proposition \ref{prop:SI1}, this means 
$C_1^2 = C_2^2 = 0$.
We also have
$$
4 = (C_1 + C_2)^2 = C_1^2 + C_2^2 + 2C_1C_2.
$$
So we obtain $C_1C_2=2$.
Thus, the curve $C=C_1+C_2$ is tangent to each other at $p$ with 
order two.
\proofend

\medskip
By using the double covering structure on some Segre surfaces over a quadric surface or a quadric cone in $\CP_3$ given in \cite{H20}, 
it is possible to see that the situation as in the two items
of Proposition \ref{prop:2mem} really happen.
For the situation in the first item, the conic with multiplicity two, for which we simply call a {\em double conic}, satisfies the following property.

\begin{proposition}\label{prop:dcp}
Irreducible double conics on a Segre surface
always form a pencil on the surface.
Further, any double conic passes a singular point of the surface.
%If a hyperplane section of a Segre surface $S$ is of the form
%$2C$ for an irreducible conic $C$, then 
%$2C$ is a 
%Assume that a hyperplane section $H|_S$ of a Segre surface $S$
%is an irreducible double conic. Then it passes a singular point of $S$.
\end{proposition}

\proof
Write the double conic as $2C$. 
Let $\mu:\tilde S\to S$ be the minimal resolution of 
all singularities of $S$ and $\tilde C$ the strict transform
of $C$ into $\tilde S$.
Then by Proposition \ref{prop:SI1}, we have $\tilde C^2 = 0$.
Hence we have $\dim |\tilde C| = 1$.
Mapping twice of members of the pencil $|\tilde C|$
by $\mu$, we obtain the required pencil of
double conics on $S$.

For the latter assertion, writing $H|_S = 2C$ with
a hyperplane $H$, we have
$$
4 = (H|_S)^2 = (2C)^2 = 4C^2.
$$
This means $C^2=1$. On the other hand, 
if the conic $C$ does not pass any singular point of $S$, then 
we have $C^2 = 0$ by Proposition \ref{prop:SI1}.
Hence the conic $C$ passes a singular point of $S$.
\proofend

\medskip
Next we show the following proposition meaning that, among the three possibilities 
in Proposition \ref{prop:2mem} for the two hyperplane sections,
the last possibility cannot occur if the point $p$ is generic.

\begin{proposition}\label{prop:vA3}
Let $S$ be any Segre surface and
 $C$ a hyperplane section of 
$S$ which is as in the second item in Proposition \ref{prop:2mem}.
Then the parameter space of the versal family of equi-singular 
displacements of $C$ in $S$
is smooth and 1-dimensional. 
\end{proposition}

\proof
The idea of the proof is similar to Theorem \ref{thm:cuspcusp}. Write $C = C_1 + C_2$ as above.
So $C_1$ and $C_2$ are irreducible conics which are tangent 
to each other at a point $p\in S_{\reg}$ with order two, and 
which have no other intersection.
As in the proof of Proposition \ref{prop:2mem},
$C_1$ and $C_2$ do not pass any singularity of $S$.
We may replace the cuspidal curve $C$ in the short exact sequence \eqref{N'1}
by the reducible curve $C_1+C_2$, 
and we still have an exact sequence
\begin{align}\label{N'2}
0 \lras \ms I_{C_1 + C_2}/\ms I_{C_1 + C_2}^2 \lras \Omega_S|_{C_1 + C_2}\lras \Omega_{C_1 + C_2}\lras 0.
\end{align}
Applying $\ms H\!om_{\ms O_{C_1 + C_2}}(\,\cdot\,,\ms O_{C_1 + C_2})$ to this sequence, 
the equi-singular normal sheaf $N'_{C_1 + C_2}$ is again defined 
as the kernel sheaf of 
the surjective homomorphism $N_{C_1 + C_2} \lras T^1_{C_1 + C_2}$,
and we have an exact sequence
\begin{align}\label{eqnn3}
0\lras N'_{C_1 + C_2} \lras N_{C_1 + C_2} \lras T^1_{C_1 + C_2} \lras 0.
\end{align}
%In order to show the assertion of the proposition,
%it is enough to show that $H^1(N'_{C_1 + C_2}) = 0$ and $H^0(N'_{C_1 + C_2})\simeq\CC$.
Let $\nu:\tilde C_1\sqcup \tilde C_2 \to C_1\cup C_2$ be the normalization 
of the curve $C_1 + C_2$.
$\tilde C_1$ and $\tilde C_2$ are isomorphic to $C_1$ and $C_2$
respectively, and the point $p=C_2\cap C_2$ is an $A_3$-singularity (i.e.\,tacnode) of $C_1 + C_2$. So 
in a neighborhood of $p$,
the curve $C_1 + C_2$ can be identified with
the curve $y^2 = x^4$ in $\CC^2$ around the origin.
As  $C_1^2 = C_2^2 = 0$ from Proposition \ref{prop:SI1}
and $C_1.\,C_2 = 2$,
we have 
\begin{align}\label{nc12}
N_{C_1 + C_2}|_{C_i}\simeq [C_1 + C_2]|_{C_i}\simeq \ms O_{C_i}(2),\quad i = 1,2.
\end{align}

Let $\ms J$ be the Jacobian ideal sheaf of 
the curve $C_1+C_2$, and write $\ms J_{C_1+C_2}$
for $\ms J|_{C_1+C_2}$.
We readily have $\ms J_p = (x^3,y)$ for the germ at $p$.
The desingularization of the tacnode may be obtained by 
blowing-up the surface $S$ twice. 
Let $\mu$ be the composition of the two blowups.
Then the restrictions of $\mu$ to 
$\tilde C_1\sqcup\tilde C_2$ gives the normalization $\nu$ of 
the curve $C_1+C_2$.
By introducing the coordinates $u= y/x$ and $v = u/x$
on the exceptional curves of the first and the second blowup respectively,
the composition $\mu$ is concretely written as
$(x,v)\longmapsto (x, x^2v)$,
while the strict transforms $\tilde C_1$ and $\tilde C_2$ are
defined by $v=1$ and $v=-1$.
Hence, if $\tilde p_1\in\tilde C_1$ and $\tilde p_2\in\tilde C_2$ denote the points on the normalization 
which correspond to the tacnode $p$,
the germs of $\mu^*\ms J$ at $\tilde p_1$ and $\tilde p_2$
are generated by two functions $x^3$ and $x^2v$.
Restricting $\mu^*\ms J$ to the strict transforms means a substitution of
$v=\pm 1$, and therefore,
over the 
two components $\tilde C_1$ and $\tilde C_2$,
we have
\begin{align}
\big(\nu^*\ms J_{C_1+C_2}\big)\big|_{\tilde C_i} &\simeq 
\big(\mu^*\ms J\big)\big|_{\tilde C_i}\notag\\
&\simeq \ms O_{\tilde C_i}(-2\tilde p_i).\label{nuJ2}
\end{align}
As in \eqref{NJ} for the case of an $A_2$-singularity,
we have $N' \simeq N\otimes_{\ms O_{C_1+C_2}} \ms J_{C_1+C_2}.$
Hence, combined with \eqref{nc12}, over each component $\tilde C_i$, 
$\nu^* N'$ is isomorphic to 
$\ms O_{\tilde C_i} (2) \otimes \ms O _{\tilde C_i}(-2) \simeq\ms O_{\tilde C_i}$.
So we have
\begin{align}
\nu^* N'_{C_1 + C_2} &\simeq \ms O_{\tilde C_1\sqcup \tilde C_2}.
\end{align}
Therefore we obtain an exact sequence
\begin{align}\label{ses:Ls}
0 \lras  N'_{C_1 + C_2} \lras \nu_*\ms O_{\tilde C_1\sqcup \tilde C_2}
\lras \CC_p \lras 0.
\end{align}
Since $R^q\nu_*\ms O_{\tilde C_1\sqcup \tilde C_2} = 0$ for $q>0$,
the Leray spectral sequence means
$H^i(\tilde C_1\sqcup \tilde C_2,\ms O_{\tilde C_1\sqcup \tilde C_2})\simeq H^i(\nu_*\ms O_{\tilde C_1\sqcup \tilde C_2})$ for any $i\ge 0$.
This readily means $H^i(\nu_*\ms O_{\tilde C_1\sqcup \tilde C_2})=0$ for any $i>0$.
From the cohomology exact sequence of \eqref{ses:Ls},
noting that the homomorphism $H^0(\nu_*\ms O_{\tilde C_1\sqcup \tilde C_2})
\lras H^0(\CC_p)$ is evidently surjective, we obtain 
$$
H^0\big(N'_{C_1 + C_2} \big)\simeq\CC \qandq H^1\big(N'_{C_1 + C_2} \big) = 0.
$$
The former implies that any first order equi-singular displacement of
the curve $C_1+C_2$ in $S$ is unobstructed, and 
the latter means that the versal family of 
such displacements of $C_1+C_2$ in $S$ is 1-dimensional.
\proofend

\medskip
One might expect that the possibility of occurrence of the double conic as in the first item
in Proposition \ref{prop:2mem} in case the point $p$ is generic can be eliminated as well.
However, as we shall see soon, this is not at all correct.
Rather, occurrence of the double conic
for a generic point $p$ will play a key role for proving
irreducibility of the cuspidal locus in the dual variety of arbitrary Segre surface.

We still assume that $p\in S_{\reg}$ is a point 
which does not belong to any line on $S$,
and $H_1$ and $H_2$ are the two hyperplanes 
which cut out from $S$ curves having non-nodal singularity 
at $p$.
%Then as an immediate implication of Propositions \ref{prop:2mem} and \ref{prop:vA3}, we have the following criteria.
%
%\begin{proposition}
%The cuspidal locus in $S^*$ is empty iff
%both of the hyperplane sections $H_1|_S$ and $H_2|_S$ are double conics for a generic point $p$.
%If exactly one of $H_1|_S$ and $H_2|_S$ is a double conic,
%then the cuspidal locus in $S^*$ is birational to $S$.
%If none of $H_1|_S$ and $H_2|_S$ is a double conic,
%then the cuspidal locus in $S^*$ is equal to $\pi_2(\bm D_S)$,
%and is birational to a double covering over $S$.
%\end{proposition}
%
%\proof
From Propositions \ref{prop:2mem} and \ref{prop:vA3}, assuming genericity to the point $p$, % and possibly exchanging $H_1$ and $H_2$, 
there are the following three possibilities for the two hyperplane sections:
\begin{itemize}
\item[(i)]
both $H_1|_S$ and $H_2|_S$ are double conics.
\item[(ii)]
exactly one of $H_1|_S$ and $H_2|_S$ has an ordinary cusp at $p$,
and the other one is a double conic,
\item[(iii)]
both $H_1|_S$ and $H_2|_S$ have ordinary cusp at $p$.
\end{itemize}
When the situation (i) occurs,
the image $\pi_2(\bm D_S)$ is not a component of the cuspidal locus
because singularity of a double curve is 
not an ordinary cusp.
Together with Proposition \ref{prop:cuspD2}, this means that 
{\em the cuspidal locus in $S^*$ is empty} in this case.
For the case (ii), the divisor $\bm D_S$ is necessarily 
reducible. It consists of two irreducible components
since the projection $\bm D_S\to S$ is of degree two.
Both components are birational to $S$ by the projection.
The image under the projection $\pi_2:I(S)\to S^*$ of the
irreducible component corresponding to the double conics 
is not a component of the cuspidal locus by the same reason to the case (i), 
while the image of the other component is
exactly the cuspidal locus in $S^*$.
It follows that {\em the cuspidal locus in $S^*$ is irreducible} in this case.
Finally, if the situation (iii) occurs, {\em whole of the image $\pi_2(\bm D_S)$ is exactly the cuspidal locus in $S^*$}. 
But at this stage, we do not know whether it is irreducible or not,
because we do not know whether the divisor $\bm D_S$ is irreducible or not.
This ambiguity will be resolved in Sections \ref{s:scl}
and \ref{s:conc}.

\subsection{Pencils of double conics on Segre surfaces}
\label{ss:pdc}
As one can expect from results in the previous section,
irreducibility of the cuspidal locus in $S^*$ is 
closely related to the existence of double conics in $S$.
In this subsection, we first provide a relation between
the pencil of double conics
and irreducibility of the cuspidal locus in the dual varieties
of Segre surfaces.
Next we determine the number of pencils of double conics
for any Segre surfaces.

\begin{proposition}\label{prop:pendc1}
Pencils of  double conics on a Segre surface $S$ are in one-to-one correspondence
with irreducible components of the divisor $\bm D_S$ 
whose generic point $(p,H)$ satisfies the property that 
$H|_S$ is a double conic.
\end{proposition}

\proof
Let $\ms P$ be a pencil of double conics on $S$.
For a generic point $p\in S_{\reg}$, there exists
precisely one member of $\ms P$ which passes $p$.
Let $H_p$ be the hyperplane which cuts out this member.
Then %as long as $p$ is not on any line on $S$,
we have $(p,H_p)\in \bm D_S$ by Corollary \ref{cor:D1}, and 
the assignment $p\mapsto (p,H_p)$ gives a birational map
from $S$ to an irreducible component of $\bm D_S$.
If $\ms P_1$ and $\ms P_2$ are distinct pencils
of double conics on $S$, then the components of $\bm D_S$ obtained this way
are distinct because $\ms P_1$ and $\ms P_2$ do not have
a common member by Proposition \ref{prop:SI1}.
Conversely, any component of the divisor $\bm D_S$ 
as in the proposition comes from a pencil of 
double conics by Proposition \ref{prop:dcp}.
Therefore, the correspondence is bijective.
\proofend

\medskip
In particular, any Segre surface $S$ has at most two pencils of 
double conics, and if $S$ has two such pencils, then 
the divisor $\bm D_S$ is reducible.
The number of pencils of double conics on each Segre surface can be determined
from the Segre symbol of the surface by the following result.

\begin{proposition}\label{prop:pendc2}
Pencils of double conics on a Segre surface $S$
are in one-to-one correspondence with 
the units 
\begin{align}\label{units}
(11), (12),(13), (14)
\end{align}
which are contained in the Segre symbol of $S$.
%
%The number of pencils of double conics on a Segre surface $S$
%is equal to the number of the following units contained in the Segre symbol of $S$: 
%\begin{align}\label{units}
%(11), (12),(13), (14).
%\end{align}
\end{proposition}

\noindent
For example, if the Segre symbol of $S$ is $[(12)11]$, then 
$S$ has a single pencil of double conics since the symbol contains one unit $(12)$ among \eqref{units}.
Similarly, if the Segre symbol of $S$ is $[(11)(11)1]$, then 
$S$ has two pencils of double conics since the symbol contains a unit $(11)$ twice.

\proof
First we explain how we can obtain a pencil of double conics
from any one of the units \eqref{units}.
The Segre symbol determines two symmetric matrices of size $5\times 5$ which correspond to two defining quadratic polynomials for the Segre surface.
An entry in the Segre symbol determines submatrices of the two full matrices.
For the unit $(12)$ for example, the two symmetric submatrices are given by
$$P_{(12)}=
\begin{pmatrix}
\aaa   & 0       & 0        \\
0      & 0       & \aaa     \\
0      & \aaa    & 1       \\
\end{pmatrix},
\quad
Q_{(12)}=\
\begin{pmatrix}
1      & 0       & 0        \\
0      & 0       & 1     \\
0      & 1       & 0       \\  
\end{pmatrix},
$$
where $\aaa$ is a complex number.
If the Segre symbol of $S$ contains the unit $(12)$,
then the two full symmetric matrices $P$ and $Q$ which define $S$
have the matrices $P_{(12)}$ and $Q_{(12)}$ as submatrices
respectively.
Therefore, the matrix $P - \aaa Q$ contains $P_{(12)} - \aaa Q_{(12)}$ as a submatrix, and we have
$$
P_{(12)} - \aaa Q_{(12)}
=
\begin{pmatrix}
0      & 0       & 0        \\
0      & 0       & 0    \\
0      & 0       & 1       \\  
\end{pmatrix}.
$$
This means that a projection from $\CP_4$ to $\CP_2$
which drops two coordinates maps $S$ to a conic in $\CP_2$.
The conic has to be irreducible since $S$ is irreducible.
By pulling back tangent lines of the conic by the last projection,
we obtain a pencil of double conics on $S$.
The same argument works %when $S$ has a pencil of double conics
when the Segre symbol contains any one of the units $(11),(13)$ or $(14)$.
(For the case $(14)$, the submatrices are equal to the full matrices.)

%
%First we show that if the Segre symbol of $S$ contains
%a unit $(14)$, then $S$ has a pencil of double conics.
%Let $\aaa$ and $\bbb$ be mutually different complex numbers,
%and define two symmetric matrices $P$ and $Q$ by 
%$$P=
%\begin{pmatrix}
%\aaa   & 0       & 0      & 0      & 0\\
%0      & 0       & 0      & 0      & \aaa\\
%0      & 0       & 0      & \aaa   & 1 \\
%0      & 0       & \aaa   & 1      & 0\\ 
%0      & \aaa    & 1      & 0      & 0
%\end{pmatrix},
%\quad
%Q=\
%\begin{pmatrix}
%1      & 0       & 0      & 0      & 0\\
%0      & 0       & 0      & 0      & 1\\
%0      & 0       & 0      & 1      & 0 \\
%0      & 0       & 1      & 0      & 0\\ 
%0      & 1       & 0      & 0      & 0
%\end{pmatrix},
%$$
%Let
%$f(X_0,\dots, X_4)$ and $g(X_0,\dots, X_4)$ be
%quadratic polynomials which correspond to $P$ and $Q$ respectively
%in a standard way.
%From the definition of Segre symbol, the Segre surface
%$S$ is defined by $F=G=0$.
%Then we have
%$$
%P - \aaa Q = 
%\begin{pmatrix}
%0      & 0       & 0      & 0      & 0\\
%0      & 0       & 0      & 0      & 0\\
%0      & 0       & 0      & 0      & 1 \\
%0      & 0       & 0      & 1      & 0\\ 
%0      & 0       & 1      & 0      & 0
%\end{pmatrix}.
%$$
%The quadratic polynomial corresponding to this symmetric matrix
%does not contain monomials involving $X_0$ and $X_1$.
%This means that the projection $(X_0,\dots, X_4)\mapsto (X_2,X_3,X_4)$ from $\CP_4$ to $\CP_2$ maps $S$
%to a conic $X_2^2 + X_3^2 + X_4^2 = 0$.

In order to provide the reverse direction, assume that a Segre surface $S$ 
has a pencil $\ms P = \{2C_t\set t\in\CP_1\}$ of double conics.
Then pulling back the conics to the minimal resolution $\tilde S$ of $S$ as in 
the proof of Proposition \ref{prop:dcp},
we obtain a pencil $\{\tilde C_t\set t\in \CP_1\}$
of rational curves on $\tilde S$
whose self-intersection numbers are zero.
The linear system 
$\{\tilde C_t + \tilde C_s\set t,s\in \CP_1\}$
induces a rational map $\tilde S\to \CP_2$ whose image is 
an irreducible conic, and this factors through $S$ via the minimal resolution.
Hence there is a rational map $S\to \CP_2$ whose image is an irreducible conic.
Since the curves $C_t + C_s$ belong to the class of hyperplane sections,
the last projection $S\to \CP_2$ is induced 
from a linear projection $\CP_4\to\CP_2$ as a restriction of the domain.
Let $f(X_0,\dots, X_4)$ and $g(X_0,\dots, X_4)$ be quadratic 
polynomials that define $S$, and $P$ and $Q$ the symmetric matrices
that correspond to $F$ and $G$ respectively.
Then from the presence of the projection to the conic, 
there has to exist constants $a$ and $b$
such that $\rank\,(aF + bG) = 3$.
Now, using the list of the symbols for all Segre surfaces
(see Table \ref{table1}),
it is not difficult to verify that this occurs only when 
the symbol includes one of $(11), (12),(13)$ and $(14)$.
\proofend

\medskip
From the proof, it follows that the component of $\bm D_S$ as in 
Proposition \ref{prop:pendc1} is contracted to a conic in $S^*$ by
the projection $\pi_2:I(S)\to S^*$.

Proposition \ref{prop:pendc2} implies the following 
\begin{corollary}\label{cor:empty}
Let $S\subset\CP_4$ be a Segre surface.
The cuspidal locus in $S^*$ is empty if and only if 
the Segre symbol of $S$ is either $[(11)(11)1]$ or $[(12)(11)]$.
\end{corollary}

\proof
These two symbols actually contain two units
among \eqref{units}. Conversely, from the list of 
symbols for all Segre surfaces as in Table \ref{table1}, 
all other symbols
contain at most one unit among \eqref{units}.
\proofend

\begin{corollary}\label{cor:bir}
If the symbol of a Segre surface $S$ is among the following list,
then the cuspidal locus in $S^*$ is an irreducible surface,
and it is birational to $S$.
\begin{align}\label{list3}
[111(11)],\,\, [12(11)],\,\, [11(12)],\,\, [1(13)],\,\,
 [(11)3],\,\, [(12)2],\,\, [(14)].
\end{align}
\end{corollary}

\proof
From the list of 
symbols for Segre surfaces (Table \ref{table1}),
these seven symbols are all ones which have exactly one unit
among \eqref{units}.
By Proposition \ref{prop:pendc1}, the divisor $\bm D_S$ for these Segre surfaces is reducible, and exactly one of the irreducible components
is from the pencil of double conics.
From the results in the previous subsection,
the image of the other irreducible component of $\bm D_S$ by the projection
$\pi_2:I(S)\to S^*$ is exactly the cuspidal locus in $S^*$.
The projection from that irreducible component of $\bm D_S$ to the cuspidal locus is birational, since the hyperplane section of $S$ determined from 
a generic point of the component has a single cusp, which means that 
the projection is generically one-to-one.

\proofend
\medskip

\begin{table}[h]
\begin{tabular}{|c|c|c|c|c|c|c|c|}
\hline
Segre symbol & 
$\Sing (S)$ & $\#$\,lines  & $\bm D_S$ & $x$ &
$y$  \\
\hline
[11111] & none & $16 + 0 + 0$ & irreducible & - & $16+0+0$  \\
\hline
[1112] & ${{\rm A}}_1$ & $8 + 4 + 0$ & irreducible  &  - & $8+ 2\cdot 4+ 0$\\
\hline
[111(11)] & $2{{\rm A}}_1$ & $0 + 8 + 0$ & reducible & - & $0 + 2\cdot 8  + 0$\\
\hline
[12(11)] & $3{{\rm A}}_1$ & $0 + 4 + 2$ & reducible  & 2 & $0 + 2\cdot 4+ 4\cdot 2$ \\
\hline
[1(11)(11)] & $4{{\rm A}}_1$ & $0 + 0 + 4$ & $\emptyset$ & 2 & $0 + 0 + 4\cdot 4$   \\
\hline
$[113]$ & ${{\rm A}}_2$ & $4 + 4 + 0$ & irreducible & - & 
$4 + 3\cdot 4 + 0$ \\
\hline
$[122]$ & $2{{\rm A}}_1$ & $ 4 + 4 + 1$ & irreducible & 2  & 
$ 4 + 2\cdot 4+4\cdot 1$\\
\hline
$[11(12)]$ & ${{\rm A}}_3$  &$0 + 4 + 0$ & reducible& - &
$0 + 4\cdot 4 + 0$\\
\hline
$[14]$ & ${{\rm A}}_3$ & $2 + 3 + 0$ & irreducible & - & 
$2 + (4+4+6)+0$\\
\hline
$[1(13)]$ & ${{\rm D}}_4$ & $0 + 2 + 0$ & reducible & -  & 
$0 + 8\cdot 2 + 0$ \\
\hline
$[(11)3]$ & $2{{\rm A}}_1+{{\rm A}}_2$ & $0 + 2 + 2$  & reducible  & 3 &$ 0 + 2\cdot 2 + 6\cdot 2$\\
\hline
$[(12)2]$ & ${{\rm A}}_1+{{\rm A}}_3$  & $0 + 2 + 1$  & reducible & 4 & $0 + 4\cdot 2 + 8\cdot 1$\\
\hline
$[(11)(12)]$  & $2{{\rm A}}_1+{{\rm A}}_3$ & $0 + 0 + 2$ & $\emptyset$ & 4  & $0 + 0 + 8\cdot 2$
\\
\hline
$[(14)]$ & ${{\rm D}}_5$ & $0 + 1 + 0$ & reducible  
& - & $0 + 16\cdot 1 + 0$\\
\hline
[23] & ${{\rm A}}_1 + {{\rm A}}_2$ & $2 + 3 + 1$ & irreducible & 3 & $2+ (2+ 2 + 3)+6\cdot 1$  \\
\hline
[5] & ${{\rm A}}_4$ & $1 + 2 + 0$ & irreducible & - & 
$1+ (8+7) + 0$\\
\hline
\end{tabular}
\bigskip
\caption{See Remark \ref{rm:table} for details.}\label{table1}
\end{table}

\begin{remark}\label{rm:table}
{\em
This is a remark for Table \ref{table1}.
For the number of lines, the first number is that of 
lines disjoint from singularities of $S$.
The second one is the number of 
lines which pass exactly one singularity of $S$.
The third one is the number of 
lines which pass two singularities of $S$.
The $x$-column is the multiplicity of the 
component $\pi_1\inv(l)$ in the divisor $\bm D_1$, where $l$ is a line passing
two singularities of $S$.
These are calculated in Section \ref{s:scl2}, and 
are independent of a choice of such a line $l$.
Finally, the $y$-column shows multiplicities of 
lines,
counted from a realization of each Segre surface as
a degeneration from smooth ones.
For instance, for the surface with symbol $[12(11)]$,
$0 + 2\cdot 4+ 4\cdot 2$ means that the multiplicities
of the 4 lines which pass one singularity is two,
and those of the two lines which pass two singularities are 4.
We are not presenting multiplicities of lines which do not
pass any singularity, because they are all one.
See the remark at the end of this article for more details
and as to why 
we present these multiplicities.
\proofend}
\end{remark}

From Corollaries \ref{cor:empty} and \ref{cor:bir},
the structure of the cuspidal locus in the dual variety of 
a Segre surface is well-understood when the surface
admits a pencil of double conics.
When the surface does not admit a pencil of double conics,
we already know that the cuspidal locus is exactly the image
of the divisor $\bm D_S$ under the projection $\pi_2:I(S)\to S^*$.
The projection from the divisor $\bm D_S$ to the cuspidal locus
is birational, by the same reason to the similar assertion 
in Corollary \ref{cor:bir}.
On the other hand, $\bm D_S$ has a structure of generically finite
double covering over $S$ by the projection $\pi_1:\bm D_S\to S$.
Thus, if $S$ is a Segre surface which does not have a pencil of double conics, then the cuspidal locus in $S^*$ is birational to 
a double covering over $S$.
In the next section, we discuss the branch divisor of 
this covering for arbitrary Segre surface $S$.

\section{Lines on Segre surfaces as branch divisors}
\label{s:scl}

We begin with the following

\begin{proposition}\label{prop:brline}
For any Segre surface $S$, every component of the branch divisor of the generically finite double covering $\bm D_S\to S$ is a line lying on $S$.
\end{proposition}

\proof
%This is an immediate consequence of Proposition \ref{prop:H}.
By Proposition \ref{prop:H}, the function $\bf H$ on $I(S)$
which defines the divisor $\bm D$ does not vanish identically on 
the fiber over a point $p\in S_{\reg}$ if $p$ is not on a line on $S$,
and moreover, the quadratic polynomial ${\bf H}(p)$ does not have a double root.
Let $(p,H_1)$ and $(p,H_2)$ be the two points
corresponding to two roots of the equation ${\bf H}(p)=0$.
Then we have $(p,H_1), (p,H_2)\in \bm D$ but these do not belong to $\bm D_1$ because ${\bf H}(p)\not\equiv 0$ as above.
Hence $(p,H_1), (p,H_2)\in \bm D_S$.
Since $H_1\neq H_2$, this means that $p$ is not a branch point of the projection 
$\bm D_S\to S$.
\proofend

\medskip
We note that the proposition does not assert that a line on $S$ is always a branch divisor.
Indeed, as mentioned in Section \ref{ss:cuspidal}, a sufficient condition for a line $l$ on $S$ to be a branch divisor is that not all 
three coefficient functions of $\bf H$ vanish identically on $l$.
If this condition is not satisfied, then $l$ is a branch divisor iff
the discriminant of the quadratic function ${\bf H}/y^m$ vanishes on $l$ identically,
where $y=0$ is a local equation of the line $l$ on $S$, and $m$ is the 
minimal vanishing order of the three functions along $l$.
In this sections, we determine whether a line on $S$
is really a branch divisor of the generically finite double covering $\bm D_S\to S$.
Conclusions will vary depending on the number of singularities of $S$ belonging to the line.

\subsection{The case $l\cap\Sing(S)=\emptyset$}
\label{ss:line1}
We begin with some property for lines on $S$ which do not pass any singularity of $S$.
The number of such lines on each Segre surface is displayed in 
Table \ref{table1}. (See Remark \ref{rm:table}.)

\begin{lemma}\label{lemma:line1}
Let $S\subset\CP_4$ be a Segre surface and 
$l$ a line lying on $S$. 
Suppose $l\cap\Sing(S)=\emptyset$. Then there exists
no hyperplane $H\subset\CP_4$ such that 
the divisor $H|_S$ includes the double line $2l$ as a sub-divisor.
Further, if a point $p$ on $l$ is not on another line on $S$, then there exists no hyperplane section 
which has a triple point at $p$.
\end{lemma}

\proof
For the first assertion, take any hyperplane $H$ containing $l$,
and write $H|_S$ in the form $l + D$ with $D$ a cubic curve.
It suffices to show that $D$ does not include the line $l$
as a component.
The indeterminacy of the projection from the line $l$ is 
eliminated
by blowing up $\CP_4$ along $l$, 
and we obtain a morphism from the blowup $\tilde{\CP}_4$
to $\CP_2$.
As $l\cap\Sing(S)=\emptyset$, the strict transform $\tilde S$ of $S$ into $\tilde{\CP}_4$
is isomorphic to $S$.
Let  $\phi:\tilde S\to \CP_2$ be the restriction of the morphism
$\tilde{\CP}_4\to\CP_2$ to $\tilde S$.
Under the identification $\tilde S\simeq S$, 
$\phi$ is naturally identified with 
the map from $S$ induced by a 2-dimensional 
subsystem of $|D|$, and this subsystem is base point free.
Since $l$ is a $(-1)$-curve from Proposition \ref{prop:SI1}, 
we have
\begin{align}\label{D2}
D^2 = (K_S\inv - l)^2 = 4-2-1=1 \qandq
D.\,l = (K_S\inv - l).\,l = 1-(-1) =2.
\end{align}
The former means that the morphism $\phi:\tilde S\to \CP_2$ is birational.
Hence the restriction $\phi|_l$ cannot be of degree-two,
and hence, from the latter of \eqref{D2}, the image $\phi(l)$ is a conic.
This means that no member of the above 2-dimensional 
subsystem of $|D|$ can contain $l$ as a component.

For the second assertion, suppose that a point $p\in l$ is not on
another line on $S$,
and let $H|_S$ be a hyperplane section which has
a triple point at $p$.
Then the hyperplane $H$ has to contain $T_pS$.
Hence we may write $H|_S = l + D$ as before, where $D$ is a possibly reducible cubic curve which has at least a double point at $p$.
If $D$ is reducible or non-reduced, then $D$ has to include a line $l'$
as a component. 
We have $l'\neq l$ from what we have already proved in this proof,
and also we have $p\not\in l'$ as we are assuming that
$p$ is not on another line on $S$.
This means that the curve $D - l'$ is a 
conic which has a double point at $p$.
So $D-l'$ consists of two lines and both of them pass $p$.
This again contradicts the assumption that 
$p$ is not on another line on $S$.
Therefore the cubic curve $D$ has to be non-reduced, irreducible,
and singular.
This means that $D$ is planer, and contradicts Proposition \ref{prop:SI2}.
Therefore no hyperplane section of $S$ can have
a triple point at $p$.
\proofend

\medskip
By using this lemma, we next show:
% that, if a line
%does not pass any singularity of $S$, then it is
%really a branch divisor of the generically finite double covering $\bm D_S\to S$,
%and moreover, it is a simple branch divisor:

\begin{proposition}\label{prop:line1}
Let $S\subset\CP_4$ be a Segre surface.
If $l$ is a line on $S$ that does not intersect $\Sing(S)$, then 
$\pi_1\inv(l)\not\subset \bm D_1$ holds, and
$l$ is a simple branch divisor of the generically finite double covering
$\bm D_S\to S$.
\end{proposition}

\proof
First we show that if $l$ is a line as in the proposition,
then for any point $p\in l$ which is not on
another line on $S$, any member of
the residual pencil 
$$
(T_pS)^*|_S - l
=\big\{D-l\set D\in (T_pS)^*|_S\big\}
$$ 
is smooth at $p$, and any two different members of the same pencil intersect each other transversely at $p$.
The former follows immediately from Lemma \ref{lemma:line1},
because a member of the residual pencil 
which is singular at $p$ gives a hyperplane section which
has a triple point at $p$.
For the latter,
let $\mu:S'\to S$ be a blowup at $p$, and $E$
the exceptional curve.
Again by Lemma \ref{lemma:line1}, the restriction homomorphism
\begin{align}\label{rst1}
r_E:H^0 \big(\mu^* H - 2 E \big) \lras H^0 \big(\ms O_E(2) \big)
\end{align}
is injective, and therefore, has a 2-dimensional image.
Let $l'$ be the strict transform of $l$ into $S'$.
Any element of the image of $r_E$ vanishes 
at the point $E\cap l'$.
So this point is a base point of the pencil on $E$
formed by the image of $r_E$.
These mean that members of the pencil $(T_pS)^*|_S - l$ are naturally in one-to-one
correspondence with points on $E$.
Therefore, any two different members of the pencil 
$(T_pS)^*|_S - l$ intersect each other transversely at $p$.

We still fix any point $p\in l$
which is not on another line on $S$.
Since $S$ is smooth at $p$ and $l$ is a line,
we can choose non-homogeneous coordinates $(x,y,z,w)$ on $\CP_4$, centered at the point $p$,
such that 
\begin{align}\label{TpS}
T_p S = \{z=w=0\}\qandq l= \{y=z=w=0\}
\end{align}
hold.
The pair $(x,y)$ can be used as coordinates on $S$
around $p$.
As in Section \ref{ss:cuspidal},
there exists a pair of holomorphic functions $F(x,y)$ and $G(x,y)$ such that 
\begin{align}\label{Seq}
S = \{(x,y,z,w)\set z = F(x,y),\,w = G(x,y)\}
\end{align}
holds around $p$.
Then  $F(x,y)=0$ and $G(x,y)=0$ are local equations around $p$ 
of generators of the pencil $(T_pS)^*|_S$.
Let $\mathfrak m=(x,y)$ be the maximal ideal at $p$.
Then we have $F,G\in \mathfrak m^2$,
and also $F(x,0)=0$ and $G(x,0)=0$ hold since $l\subset S$.
So we can write $F = yf$ and $G = yg$ for some $f,g\in \mathfrak m$.
From the property on members of the residual pencil 
$(T_pS)^*|_S - l$ obtained in the first part of 
the present proof, we may suppose that
the curve $\{f=0\}$ intersects $l$ transversely at $p$.
This condition implies $f_x(0,0)\neq 0$.
We consider the map 
$$
(x,y)\longmapsto (\tilde x, \tilde y) := \big( f(x,y), y\big)
$$
defined in a neighborhood of $p$ in $S$.
As $f_x(0,0)\neq 0$, the Jacobian of this map 
does not vanish at the origin, and hence we can use
$(\tilde x,\tilde y)$ as holomorphic coordinates
in a neighborhood of $p$.
For simplicity of presentation, we redefine $x$
to be the above $\tilde x$, so that we have $F= xy$,
and we still write $G = y g$, $g\in\mathfrak m$.
Then by replacing $G$ by $G - c F$ for a constant $c$
to eliminate the term $xy$ from $G=yg$,
we may suppose that the leading term of $G$ in
a Taylor expansion at $p$ is of the form $ay^2$
for some constant $a$. 
We have $a\neq 0$ since otherwise
the zero of $G$ would have a triple point at $p$,
which contradicts Lemma \ref{lemma:line1}.
Thus, an expansion of $g\,(=G/y)$ at $p$ is 
of the form 
\begin{align}\label{expansion1}
ay + h(x,y)\quad{\text{with $h\in \mathfrak m^2$
and $a\neq 0$.}}
\end{align}

Now recall from \eqref{Hess15} that a defining equation of 
the divisor $\bm D = \bm D_1 + \bm D_S$ on $I(S)$ is given by
\begin{align}\label{Hess8}
 \Hess(F)\lmd^2
+ \big(F_{xx}G_{yy} + G_{xx}F_{yy} - 2F_{xy}G_{xy}\big)\lmd\mu
+ 
\Hess(G)\mu^2
 = 0.
\end{align}
Thanks to the above change of coordinates, we  have
$\Hess (F) = -1\neq 0$, and this means that 
$\pi_1\inv(l)\not\subset \bm D_1$.
We recall from \eqref{disc} that 
the discriminant of \eqref{Hess8} is given by
\begin{align}\label{disc2}
 \big(F_{xx} G_{yy} - G_{xx} F_{yy} \big)^2
 + 4 \big(F_{xx}G_{xy} - F_{xy}G_{xx}\big)
 \big(F_{yy}G_{xy} - F_{xy}G_{yy}\big).
\end{align}
In the present situation, it is easy to see that 
this is written as 
\begin{align}\label{disc20}
4G_{xx}G_{yy}=
4yg_{xx}(2g_y + yg_{yy}).
\end{align}
This means that the line $l=\{y=0\}$ is included in the zero divisor
of the discriminant \eqref{disc2} with multiplicity at least one.
If this multiplicity is strictly greater than one,
from \eqref{disc20}, we obtain a divisibility $y \set g_{xx}$ or $y\set g_y$.
The former happens iff the Taylor expansion of 
$g$ at the origin is of the form
$$
b_1 + b_2x + b_3y+ y\tilde h(x,y),
$$
with $\tilde h\in\mathfrak m$ and $b_1,b_2,b_3\in\CC$.
Comparing this with \eqref{expansion1},
we obtain $b_1 = b_2 = 0$ and $b_3 = a$.
This means $g = a y +y\tilde h= y(a + \tilde h)$.
Hence $y\set g$.
Therefore $y^2 \set G$ holds, which implies that 
the hyperplane section $\{w=0\}|_S$ contains the double line 
$2l$.
This contradicts Lemma \ref{lemma:line1}.
Hence we have $y \not|\, g_{xx}$.
Furthermore, again from \eqref{expansion1}, we have $g_y(0,0)\neq 0$.
Hence we have $y\not|\, g_y$.
Thus we have shown that 
the discriminant \eqref{disc20} is not divisible by 
$y^2$.
Hence the line $l$ is included
in the zero divisor of the discriminant 
with multiplicity precisely one.
This means that $l$ is a simple branch divisor of $\bm D_S\to S$.
\proofend

\medskip
We postpone an immediate implication from Proposition \ref{prop:line1} until Section \ref{s:conc}.

\subsection{The case where $l\cap\Sing(S)$ is a single point}
\label{ss:line2}
In this subsection, we consider lines on Segre surfaces to which only one singular point
of the surface belongs.
Again, the number of such lines on each Segre surface is displayed in 
Table \ref{table1}. (See Remark \ref{rm:table}.)
Each of these lines are not Cartier divisors on the surfaces.
% we show that these lines are 
%also a branch divisor of the covering $\bm D_S\to S$.
%The result will be quite different from the one obtained in the last subsection.
We begin with the following lemma which is analogous to 
Lemma \ref{lemma:line1}.

\begin{lemma}\label{lemma:line2}
Let $S\subset\CP_4$ be a Segre surface and 
$l$ a line  on $S$. 
If $l\cap\Sing(S)$ consists of one point, then there exists
a unique hyperplane $H$ such that 
the divisor $H|_S$ is of the form $2l + l' + l''$,
where $l'$ and $l''$ are possibly identical lines which are different from $l$.
Further, if a point $p\in l\cap S_{\reg}$ does not belong to any other line
on $S$, 
then there exists no hyperplane section 
which has a triple point at $p$.
\end{lemma}

\proof
Let $\mu:\tilde S\to S$ be the minimal resolution of 
the unique singularity of $S$ lying on $l$,
$E$ the exceptional divisor,
and $\tilde l$ the strict transform of $l$ into $\tilde S$.
Of course, $E$ is not necessarily irreducible.
Let $Z$ be the fundamental cycle in the exceptional divisor $E$, so that 
if $\mathfrak m$ is the maximal ideal at the singularity,
for a generic element $f\in\mathfrak m$, the function 
$\mu^*f$ vanishes on $E$
with the same multiplicity
as $Z$ on each component of $E$.
Then as in the same way to showing the property 
$\mu^*\mathfrak m\simeq \ms O(-Z)$ \cite[III (3.8) Proposition]{BHPV}, we have
\begin{align*}
\mu^*\left(K_S\inv\otimes\ms O_S(-l)\right)
\simeq
K_{\tilde S}\inv\otimes \ms O_{\tilde S}(-\tilde l - Z).
\end{align*}
Let $\tilde{\CP}_4\to \CP_4$ be a blowup along $l$,
and $S'\subset\tilde{\CP}_4$ the strict transform of 
$S$.
This time the projection $S'\to S$ is not an isomorphism,
but the minimal resolution $\mu$ factors through 
the morphism $S'\to S$.
This means that the system 
$\big| K_{\tilde S}\inv-\tilde l - Z\big|$
does not have a base point.
Using $\tilde l^2 = -1$
from Proposition \ref{prop:SI1}
as well as the properties $Z^2 = -2$ and $Z.\,\tilde l = 1$,
we calculate
\begin{align*}
\big(K_{\tilde S}\inv-\tilde l - Z\big)^2 &= 
(K_S\inv)^2 + \tilde l^2 + Z^2 - 2
K_{\tilde S}\inv.\, \tilde l + 2\tilde l.\,Z
-2 K_{\tilde S}\inv.\,Z \\
&= 4 + (-1) + (-2) -2.\,1 + 2.\,1 - 2.\,0=1,\\
\big(K_{\tilde S}\inv-\tilde l - Z\big).\, \tilde l 
&= 1 - (-1)  - 1 = 1.
\end{align*}
If $\phi:\tilde S\to \CP_2$ means 
the morphism induced by the system $\big| K_{\tilde S}\inv-\tilde l - Z\big|$, 
these imply that $\phi$ is birational and 
the image $\phi(\tilde l)$ is a line.
Thus, the divisor $\phi\inv(\phi(l))$ is the
 unique member of
$\big| K_{\tilde S}\inv-\tilde l - Z\big|$
which has $\tilde l$ as a component.
The multiplicity of $\tilde l$ in this member is precisely one since 
$\phi$ is birational.
Adding $\tilde l +Z$ to $\phi\inv(\phi(l))$ and 
projecting it to $S$ by $\mu$,
we obtain the unique hyperplane section of $S$
which includes the line $l$ with multiplicity precisely two.
Write this hyperplane section as $2l + C$, where $C$ is a conic.
This conic has to be reducible since otherwise $C$ will move
in $S$ from Proposition \ref{prop:SI1},
which contradicts the uniqueness just proved.
Therefore $C$ consists of two lines,
and the unique hyperplane section is of the form
$2l + l' + l''$ with $l'\neq l$ and $l''\neq l$
as in the lemma.

For the latter assertion of the lemma,
let $p\in l\cap S_{\reg}$ be as in the lemma,
and suppose that $H|_S$ is a hyperplane section
which has a triple point at $p$.
The hyperplane section $2l + l' + l''$ has a triple point,
but the triple point is on a line different from $l$.
So we have $H|_S\neq 2l + l' + l''$.
Hence, from the above uniqueness, $H|_S$ is of the form
$l + D$, where $D$ is a cubic curve which does not contain $l$
as a component. $D$ has at least a double point at $p$.
If $D$ is reducible, it contains a line, which does not 
pass $p$. But then the residual conic $D-l$ cannot have
a double point at $p$. Hence $D$ is irreducible singular cubic curve.
Such a curve is planer.
This contradicts Proposition \ref{prop:SI2}.
Therefore, no hyperplane section has a triple point at $p$.
\proofend

\medskip By using this lemma, we show an analogous result to 
Proposition \ref{prop:line1}.

\begin{proposition}\label{prop:line2}
Let $S\subset\CP_4$ be a Segre surface and 
$l$ a line on $S$.
If $l\cap \Sing (S)$ consists of one point, then 
$\pi_1\inv(l)\not\subset \bm D_1$ holds,
and
$l$ is a branch divisor of the generically finite double covering
$\bm D_S\to S$ with multiplicity at least two.
\end{proposition}

\proof
The idea of the proof is similar to Proposition \ref{prop:line1}.
First we show that if the line $l$ is as in the proposition,
then for any  point $p\in l\cap S_{\reg}$ which is not on
another line on $S$, we can choose two generators of 
the residual pencil $(T_pS)^*|_S - l$ in such a way
that one of them has the line $l$ as the unique component
which passes $p$, while the other one has a component
which intersects $l$ transversely at $p$
as  the unique component which passes $p$.

First let $H_1|_S$ be the unique hyperplane section
which includes the double line $2l$ as in 
the first half of Lemma \ref{lemma:line2}.
By the assumption that $p$ does not belong to another
line on $S$ and the latter half of Lemma \ref{lemma:line2},
there does not exist a hyperplane section of $S$ 
which has a triple point at $p$.
Hence the residual curve $H_1|_S - 2l$ does not pass $p$.
Let $H_2$ be any element of the pencil $(T_pS)^*$
which is different from $H_1$.
Then since the multiplicity of $H_2|_S$ at $p$ is 
precisely two by Lemma \ref{lemma:line2},
the curve $H_2|_S- l$ has a unique component
which passes $p$, and it is smooth at $p$.
If this component were tangent to $l$ at $p$,
by subtracting a constant multiple of 
a defining equation of $H_1$ from that of $H_2$,
there would exist a hyperplane $H\in (T_pS)^*$
for which $H|_S$ has a triple point at $p$.
This contradicts Lemma \ref{lemma:line2}.
Hence the component always intersects $l$ transversely at $p$.
Thus we have obtained generators of 
the pencil $(T_pS)^*|_S - l$ which satisfy
the above properties.

As in the proof of Proposition \ref{prop:line1}, 
we can choose non-homogeneous coordinates $(x,y,z,w)$
centered at $p$, such that 
the properties \eqref{TpS} hold.
Let $F(x,y)$ and $G(x,y)$ be a pair of holomorphic functions
such that \eqref{Seq} holds.
We may suppose that the hyperplane $\{z=0\}$ is exactly the 
one whose section contains the double line $2l$.
Then we have the divisibility $y^2 \set F$ and $y \set G$.
Further, if we write $G = yg(x,y)$, then 
$g$ belongs to the maximal ideal $\mathfrak m$ at $p$.
Moreover, we have $g_x(0,0)\neq 0$ from the above 
transversality for the component different from $l$.
Hence in the same way to the proof of Proposition \ref{prop:line1},  we may use $(g(x,y), y)$ as coordinates on $S$ around $p$
instead of $(x,y)$.
Again we redefine $x$ as $g(x,y)$, so that 
$G= xy$.
We remain to write $F = y^2 f$.
Then we have $f\not\in\mathfrak m$
since all components of $\{z=0\}\cap S$ passing $p$ are $2l$.
From $G=xy$, we have $\Hess(G) = -1\neq 0$,
and hence again from \eqref{Hess8}, we have $\pi_1\inv(l)
\not\subset \bm D_1$.
Moreover, this time from $F = y^2 f$, we obtain that
the discriminant \eqref{disc2} is given by 
\begin{align}\label{disc3}
4F_{xx} F_{yy} = 4y^2 f_{xx} (2f + 4y f_y + y^2 f_{yy}).
\end{align}
Since this is divisible by $y^2$, we obtain that 
the generically finite double covering $\bm D_S\to S$ has 
the line $l$ as a branch divisor
and the multiplicity of $l$ as a branch divisor is at least two.
\proofend

\medskip
We remark that since the function $f_{xx}$ in
\eqref{disc3} may be divisible by $y$ or its power,
we cannot conclude that the multiplicity of the line $l$
as a branch divisor is exactly two. 
Note that $2f + 4y f_y + y^2 f_{yy}$ is not divisible by $y$ since $f$ is not.
Since quadric polynomials which define a Segre surface are concretely obtained, it might be possible to
obtain the functions $F$ and $G$ in explicit forms
in the present circumstance,
and by using them to determine the multiplicity 
of the line as a branch divisor, for each Segre surface
and each line of the present kind.
In Section \ref{s:scl2}, we will do that sort of calculations 
for lines on Segre surfaces which pass two singularities of 
the surfaces.

\subsection{The case where $l\cap\Sing(S)$ consists of two points}
\label{ss:line3}
%Let $S$ be a Segre surface which does not
%have a pencil of double conics, or equivalently, a Segre surface
%which has a line that does not pass any singularity.
%The divisor $\bm D_S$ is irreducible from Proposition \ref{prop:irr2},
%and a line on $S$ is a branch divisor of $\bm D_S\to S$
%if the line passes at most one singular point of $S$
%from Propositions \ref{prop:line1} and \ref{prop:line2}.
%Among these Segre surfaces, those ones which have a line
%passing two singular points of the surfaces are 
%$[221]$ and $[32]$.
Let $l$ be a line on a Segre surface $S$ which passes two singular
points of $S$.
In this subsection, we investigate the structure of hyperplane sections
of $S$ which belong to the pencil $(T_pS)^*|_S$ for some point $p\in l\cap S_{\reg}$.
Note that since $S$ is quartic, no three singularities
of $S$ belong to the same line lying on $S$.
Note also that the line joining two singularities 
of $S$ does not necessarily lie on $S$.
We begin with the following property which is characteristic to
the present kind of lines.

\begin{lemma}\label{lemma:line3}
Let $S\subset\CP_4$ be a Segre surface and 
$l$ a line on $S$. 
If $l\cap\Sing(S)$ consists of two points, then 
the tangent plane $T_pS$ is independent of 
a choice of a point $p\in l\cap S_{\reg}$.
In particular, for any such a point $p$ 
and any hyperplane $H\in (T_pS)^*$,
the section $H|_S$ includes the double line $2l$ as a sub-divisor.
\end{lemma}

\proof
Let $p_1$ and $p_2$ be the two singularities of $S$ on $l$,
and $\mu:\tilde S\to S$ the minimal resolution of $p_1$
and $p_2$. 
Let $\tilde l$ be the strict transform of $l$ into $\tilde S$,
and $Z_1$ and $Z_2$ the fundamental cycles over $p_1$
and $p_2$ respectively.
Then in the same way to the proof of Lemma \ref{lemma:line2}, we have
\begin{align*}
\mu^*\left(K_S\inv\otimes\ms O_S(-l)\right)
\simeq
K_{\tilde S}\inv\otimes \ms O_{\tilde S}(-\tilde l - Z_1 - Z_2),
\end{align*}
and the system 
$\big| K_{\tilde S}\inv-\tilde l - Z_1 - Z_2\big|$
does not have a base point.
Similarly, we also have 
$\big( K_{\tilde S}\inv-\tilde l - Z_1 - Z_2\big)^2 = 1$,
but this time we have
$$
\big( K_{\tilde S}\inv-\tilde l - Z_1 - Z_2\big).\,
\tilde l = 1 - (-1)  - 1 - 1 = 0.
$$
If $\phi:\tilde S\to \CP_2$ means 
the morphism induced by the 2-dimensional system $\big| K_{\tilde S}\inv-\tilde l - Z_1 - Z_2\big|$, 
these imply that $\phi$ is a birational morphism and 
it contracts the curve $\tilde l$ to a point.
Hence the system 
$\big| K_{\tilde S}\inv- 2\tilde l - Z_1 - Z_2\big|$
(note that $\tilde l$ is further subtracted)
is a pencil.
Adding $2\tilde l + Z_1 + Z_2$ to the members of this pencil
and projecting them to $S$ by $\mu$,
we obtain a pencil of hyperplane sections
which contains the double line $2l$ as fixed components.
Clearly, this pencil is equal to the pencil 
$(T_pS)^*|_S$.
Since this holds for any $p\in l\cap S_{\reg}$,
we obtain that the pencil 
$(T_pS)^*|_S$ is independent of a choice of 
the point $p\in l\cap S_{\reg}$.
Hence the tangent plane $T_pS$ is independent of 
a choice of a  point $p\cap S_{\reg}$.
\proofend

\medskip
It is possible to obtain more detailed information
about the structure of the hyperplane sections
by $H\in (T_pS)^*$,
in the situation of Lemma \ref{lemma:line3} as follows.
We notice that from Table \ref{table1},
%\cite[Section 8.6]{Dol}, 
when two singularities of $S$ belong to
a line, the types of them are $\{A_1,A_1\},
\{A_1,A_2\}$ or $\{A_1,A_3\}$.

\begin{lemma}\label{lemma:line4}
Let $S$ and $l$ be as in Lemma \ref{lemma:line3},
and $p\in l\cap S_{\reg}$ a  point  such that 
no another line on $S$ passes $p$.
Then there exists a unique hyperplane
section of $S$ whose multiplicity at $p$ is at least three.
Further, we have:
\begin{itemize}
\item If the singularities on $l$ are $\{A_1,A_1\}$,
then such a section is  of the form $2l + C$, where $C$ is an irreducible conic that intersects $l$ transversely at $p$.
\item 
If the singularities on $l$ are $\{A_1,A_2\}$,
then such a section is  of the form
$3l + l'$ where $l'$ is a line different from $l$.
\item 
If the singularities on $l$ are $\{A_1,A_3\}$,
then such a section is of the form $4l$.
\end{itemize}
\end{lemma}

\proof
Because not every member of the pencil
$(T_pS)^*|_S$ has a triple point at $p$
by Lemma \ref{lemma:dbl},
there exists at most one member of the pencil
which has a triple point at $p$.
Let $\mu:S'\to S$ be blowup at $p$,  $E$
the exceptional curve, and $l'$ the strict transform of $l$.
By Lemma \ref{lemma:line3},
any member of the pencil $(T_pS)^*|_S$ 
contains the double line $2l$.
%Let $[l]$ be the point on $E$ which corresponds to
%the tangent direction of $l$ at $p$.
So the strict transform to $S'$ of any member
of the pencil $(T_pS)^*|_S$ has
at least a double point at $E\cap l$.
This implies that 
the image of the restriction homomorphism 
$$
r_E:H^0 \big(\mu^* H - 2 E \big) \lras H^0 \big(\ms O_E(2) \big)
$$
is at most 1-dimensional.
Further the kernel of $r_E$ is at most 1-dimensional 
from Lemma \ref{lemma:dbl}.
Hence $r_E$ has exactly 1-dimensional kernel.
This means that the pencil $(T_pS)^*|_S$
has a unique member whose multiplicity at $p$ is at least three.

Next in order to show the latter half of the proposition,
we change the meaning of the above notations and
let $\mu:\tilde S\to S, \tilde l, Z_1, Z_2$ have
the same meaning as in the proof of 
Lemma \ref{lemma:line3}. 
Further, for brevity, we put $L:= K_{\tilde S}\inv-\tilde l - Z_1 - Z_2$, and let  $\phi:\tilde S\to \CP_2$ again mean the morphism induced by the
2-dimensional linear system $|L|$.
Recall that $\phi$ is birational and contracts $\tilde l$ to a point.

First suppose that both singularities on $l$ are $A_1$-points as in the first item in the lemma. 
Then each of the two divisors $Z_1$ and $Z_2$ is just  
a single $(-2)$-curve respectively.
%As in the proof of Lemma \ref{lemma:line3},
%the image $\phi(\tilde l)$ is a point.
We readily have
$$
L.\,Z_1 = L.\,Z_2 = 1.
$$
Therefore both $\phi(Z_1)$ and $\phi(Z_2)$ are lines.
By the factorization theorem of birational morphism between smooth
surfaces,
the morphism $\phi$ is a composition of blowups.
We have $\tilde l^2 = -1$ from Proposition \ref{prop:SI1}.
Further, since $\tilde S$ is a weak del Pezzo surface, there does not exist a smooth rational curve on $\tilde S$ whose self-intersection number
is less than $(-2)$.
From these, we can readily show that no curve intersecting $\tilde l$
is contracted by $\phi$.
%
%In particular, the preimage $\phi\inv(\phi(\tilde l))$
%is a chain of smooth rational curves,
%and by Hodge index theorem,
%the intersection matrix of this chain is negative definite.
%Suppose that the chain $\phi\inv(\phi(\tilde l))$ has
%an irreducible component other than $\tilde l$,
%and let $\CCC$ be a component which intersects $\tilde l$.
%Then we have $K_{\tilde S}\inv.\,\CCC\le 1$ from adjunction formula,
%and obviously we have $\tilde l.\,\CCC\ge 1,\,
%Z_1.\,\CCC\ge 0$ and $Z_2.\,\CCC\ge 0$.
%From these we obtain $L.\,\CCC\le 0$ and equality holds
%iff $\CCC^2 = -1,\, \tilde l.\,\CCC = 1$
%and $Z_1.\,\CCC=Z_2.\,\CCC = 0$ hold.
%Since $|L|$ is base point free, we have $L.\,\CCC = 0$.
%Therefore $\CCC$ has to be a $(-1)$-curve
%which intersects $\tilde l$ transversely at a unique point.
%As $\tilde l$ is also a $(-1)$-curve,
%this contradicts that the intersection matrix of 
%the exceptional divisor $\phi\inv(\phi(\tilde l))$ is negative definite.
%Thus the divisor $\phi\inv(\phi(\tilde l))$ has no component
%other than $\tilde l$.
Hence, in a neighborhood of $\tilde l$,
$\phi$ is exactly the blow down of the $(-1)$-curve $\tilde l$.
This implies that for any line $\ell\subset \CP_2$
through the point $\phi(\tilde l)$, the preimage
$\phi\inv(\ell)$ contains the curve $\tilde l$
by multiplicity precisely one.
Hence no member of the system $|L|$
contains $\tilde l$ by multiplicity strictly greater than one.
Therefore 
no member of the system $|K_S\inv|$
contains $\tilde l$ by multiplicity strictly greater than two.
So we may write the unique member
of the pencil $(T_pS)^*|_S$
which has a triple point at the point $p$
as $2l + C$, where $C$ is a conic which goes through $p$.
Because we are assuming non-existence of a line
through $p$, the conic $C$ cannot be reducible.
If $C$ is tangent to $l$ at the point $p$,
then the curve $l + C$ would lie on the 2-plane
which is spanned by $C$, and this contradicts
Proposition \ref{prop:SI2}.
Hence the conic $C$ intersects $l$ transversely at a
unique point. This proves the first item in the proposition.

Next suppose that 
one of the two singularities on $l$ is an $A_2$-point
as in the second item of the proposition.
We may suppose that the fundamental cycle $Z_1$ is over
the $A_2$-point, and we write it as
$Z_1 = E_1 + E_2$, where $E_2$ is the component
which intersects $\tilde l$.
The cycle $Z_2$ is a single $(-2)$-curve.
We readily have
$$
L.\,\tilde l = L.\,E_2 = 0,\quad L.\,E_1 = L.\, Z_2 = 1.
$$
These mean that not only $\tilde l$ but also 
$E_2$ are contracted to a point by
the birational morphism $\phi:\tilde S\to\CP_2$,
and $E_1$ and $Z_2$ are mapped to lines by $\phi$.
%By birationality of $\phi$, we have $\phi(E_1)\neq\phi(Z_2)$.
Further we have $\phi(E_1)\neq\phi(Z_2)$ because
$\phi$ is birational.
Furthermore, in a similar way to the case of $\{A_1,A_1\}$,
using $\tilde l^2 = -1$ and $E_2^2 = -2$, 
we can show that in a neighborhood of the chain
$\tilde l + E_2$, the birational morphism $\phi$ contracts
exactly these two curves.
From the above self-intersection numbers, 
$\phi$  blows down $\tilde l$ first
and next $E_2$.
This implies that the member 
$\phi\inv\big(\phi(Z_2)\big)\in |L|$
contains $\tilde l$ by multiplicity precisely two.
Adding $\tilde l + Z_1 + Z_2$ to it and taking the image 
by $\mu:\tilde S\to S$, we obtain a hyperplane section
of $S$ which is of the form $3l + l'$,
where $l'$ is a line different from $l$.
Thus we have obtained the second item in the proposition.

Finally, suppose that one of the two singularities on $l$ is an $A_3$-point, and assume that the chain $Z_1 = E_1  + E_2 + E_3$
is the fundamental cycle over this singularity, with the component $E_3$ (resp.\,$E_2$)
intersecting $\tilde l$ (resp.\,$E_3$).
Then this time we have
$$
L.\,\tilde l = L.\,E_2 =  L.\,E_3= 0,\quad L.\,E_1 = L.\, Z_2 = 1.
$$
These mean that the chain $\tilde l + E_3+E_2$ is contracted to a point by $\phi:\tilde S\to\CP_2$,
and $E_1$ and $Z_2$ are mapped to lines by $\phi$.
Further these lines are mutually different by the same reason to
the $A_2$-case.
Since $\tilde l^2 = -1$ and $E_3^2 = E_2^2 = -2$,
the birational morphism $\phi$ blows down
components of the chain $\tilde l + E_3+E_2$ in this order.
It follows that the member $\phi\inv(\phi(Z_2))\in |L|$
contains $\tilde l$ by multiplicity precisely three.
Adding $\tilde l + Z_1 + Z_2$ to it and taking the image 
by the minimal resolution $\mu:\tilde S\to S$, we obtain a hyperplane section
of $S$ which is of the form $4l$.
Thus, we obtain the third item in the proposition.
\proofend

\medskip
Note that in the proof for the second and third items,
absence of a line through the point $p$ is not used.

By using Lemma \ref{lemma:line3}, we prove the following result
about the divisor $\bm D_1$.
%At this stage we can give a conclusion about 
%whether a line on a Segre surface is a branch divisor
%of the generically finite double covering $\bm D_S\to S$
%only for lines in the first item in
%Lemma \ref{lemma:line4}:

\begin{proposition}\label{prop:line6}
Let $S\subset\CP_4$ be a Segre surface and 
$l$ a line on $S$.
If $l\cap \Sing (S)$ consists of two points,
then $\pi_1\inv(l)$ is a component of the divisor $\bm D_1$, 
and its multiplicity is at least two.
\end{proposition}

\proof
Again an idea of the proof is similar to Proposition \ref{prop:line1}.
Pick any point $p\in l\cap S_{\reg}$.
We choose homogeneous coordinates $(x,y,z,w)$ on $\CP_3$
centered at $p$, such that the conditions \eqref{TpS} hold,
and let $F(x,y)$ and $G(x,y)$ be holomorphic functions 
around $p$ such that \eqref{Seq} holds.
By Lemma \ref{lemma:line3}, we can write $F = y^2 f$
and $G = y^2 g$ for some holomorphic functions $f$ and $g$
around $p$.
Then by elementary calculations, we have
$$
y^2\set \Hess(F),\quad
y^2\set \Hess(G) \qandq
y^2\set F_{xx}G_{yy} + G_{xx}F_{yy} - 2F_{xy}G_{xy}.
$$
Thus all coefficient functions of the defining function $\bf H$
in \eqref{Hess15} of the divisor $\bm D$ on $I(S)$
are divisible by $y^2$, and this implies that 
$\pi\inv(l)$ is contained in the divisor $\bm D_1$
by multiplicity at least two.
\proofend

\medskip
We remark that when both of the two singularities on 
the line $l$ are $A_1$-points as in the first item
in Lemma \ref{lemma:line4},
by making use of the unique hyperplane section which
has a triple point at $p\in l\cap S_{\reg}$,
it is possible to show that the multiplicity of 
the component $\pi_1\inv(l)$ in $\bm D_1$ is precisely two,
in a similar way to Propositions \ref{prop:line1} and
\ref{prop:line2}.
But if one of the two singularities on 
the line $l$ is not an $A_1$-point, Lemma \ref{lemma:line4}
does not give an element of the pencil $(T_pS)^*|_S$
whose equation can be assumed to be in a sufficiently informative form (like $xy$) as in 
the cases of Propositions \ref{prop:line1} and \ref{prop:line2},
and by this reason we cannot determine the precise multiplicity
of the component $\pi_1\inv(l)$ in $\bm D_1$.
By a similar reason, it seems more difficult to determine whether such a line $l$
is really a branch divisor of the generically finite double
covering  $\bm D_S\to S$.
In Section \ref{s:scl2}, we determine the precise multiplicity of $\pi\inv(l)$ by using explicit defining polynomials of Segre surfaces,
and show that the lines in question are  {\em always not} 
branch divisors of the double covering.

\section{Conclusions, and singularity of the cuspidal locus }
\label{s:conc}
\subsection{Irreducibility of the cuspidal locus}
\label{ss:conc}
Let $S$ be a Segre surface.
From Corollaries \ref{cor:empty} and \ref{cor:bir}, we already know that 
the cuspidal locus in $S^*$ is empty if the Segre symbol of 
$S$ is either $[(11)(11)1]$ or $[(11)(12)]$,
and that the cuspidal locus is irreducible and birational to $S$
if the symbol contains exactly one of the four units
$(11),(12),(13)$ and $(14)$. 
For the remaining Segre surfaces, it is already easy to show the following conclusion.

\begin{proposition}\label{prop:irr2}
If the Segre symbol of a Segre surface $S$ is among the following seven ones, 
then the divisor $\bm D_S$ and the cuspidal locus
in $S^*$ are irreducible, and they are mutually birational.
In particular, the cuspidal locus is birational to a double covering over $S$.
\begin{align}\label{sbd}
[11111],\,\,[1112],\,\,[113],\,\,[122],\,\,[14],\,\,[23],\,\,[5].
\end{align}
%
%Suppose that $S$ is a Segre surface which has a line not through
%any singularity of $S$, so that whole of the divisor $\bm D_S$ 
%is irreducible. Then the branch divisor of the generically double covering $\bm D_S\to S$ consists of lines not through two singular points of $S$.
\end{proposition}

\proof
The seven symbols \eqref{sbd} are characterized by the property that 
they do not contain any units among $(11),(12),(13)$ and $(14)$.
By Proposition \ref{prop:pendc2}, this condition is equivalent to 
absence of a pencil of double conics on the Segre surface.
From Proposition \ref{prop:cuspD2} and the results in Section \ref{ss:hsf}, absence of a pencil of double conics means that 
whole of the divisor $\bm D_S$ is mapped precisely 
to the cuspidal locus in $S^*$.
Further, as presented in Table \ref{table1}, the Segre surfaces having
these seven symbols are
exactly the ones which have a line that does not pass any singularity 
of the surfaces.
By Proposition \ref{prop:line1}, these lines are simple branch divisors
of the generically finite double covering $\bm D_S\to S$. 
Therefore the divisor $\bm D_S$ is irreducible for these seven kinds
of Segre surfaces.

It remains to show that the projection from $\bm D_S$ to the cuspidal locus
in $S^*$ is birational.
But this is immediate from the fact that for a generic point $(p,H)$ of the divisor $\bm D_S$, the hyperplane section $H|_S$ has a single ordinary cusp as its all singularity.
\proofend

\medskip
Irreducibility of the divisor $\bm D_S$ for these Segre surfaces are 
reflected in Table \ref{table1}.

Combined with Corollaries \ref{cor:empty} and \ref{cor:bir},
Proposition \ref{prop:irr2} gives the following result.

\begin{corollary}\label{cor:irr}
The cuspidal locus in the dual variety of any Segre surface is 
either empty or an irreducible surface.
\end{corollary}

For the structure of the cuspidal locus beyond irreducibility
for Segre surfaces in Proposition \ref{prop:irr2},
the cuspidal locus is birational to the divisor $\bm D_S$ in $I(S)$,
and $\bm D_S$ has a structure of generically finite double 
covering over $S$.
From Proposition \ref{prop:brline}, 
any component of the branch divisor of the double covering is 
a line on $S$, and from Propositions \ref{prop:line1} and \ref{prop:line2}, any line on $S$ which passes at most one singularity of $S$ is actually a branch divisor of the covering.
So what is still missing is whether lines through 
two singular points of $S$ are components
of the branch divisor of the covering $\bm D_S\to S$.

If $S$ is a Segre surface whose symbol is among \eqref{sbd}, 
and $S$ has a ling passing two singular points of the surface,
then from Table \ref{table1}, the symbol of $S$ is either $[221]$ or $[32]$.
In Appendix, we will show by explicit calculations using
defining equations of Segre surfaces, that 
the lines through two singular points of the surfaces are
{\em always not} a component of a branch divisor of $\bm D_S\to S$.
This implies the following 

\begin{proposition}\label{prop:br3}
Let $S$ be a Segre surface whose symbol is among \eqref{sbd}.
Then the branch divisor of the generically finite double covering
$\bm D_S\to S$ consists of lines which pass at most one singular point of $S$.
\end{proposition}

\medskip

Next, by using the results obtained so far, we explicitly construct
a surface which is birational to  the cuspidal locus in $S^*$ for 
a smooth Segre surface $S$.

A smooth Segre surface is nothing but a smooth del Pezzo surface 
of degree four.
It is well-known that any smooth del Pezzo surface of degree four
can be realized as 
a 5 points blown up of $\CP_2$, where 5 points are in a general position.
%and other Segre surfaces (with singularities) are obtained 
%by specializing these 5 points in suitable ways.
In the following, instead of using this realization,
we make use of a realization of smooth Segre surfaces as 4 points
blown up of $\qdr$. This seems more economical than
to work under the above realization. 

Choose 4 points on $\qdr$ such that no two points are
on the same $(1,0)$-curve nor on the same $(0,1)$-curve,
and such that there exists no $(1,1)$-curve to which all the 4 points belong.
Let $S\to \qdr$ be a blowing up at these 4 points.
Then the anti-canonical system $|K_S\inv|$ gives a projective
embedding $S\subset \CP_4$ which realizes $S$ as a smooth
complete intersection of two quadrics.
Thus, $S$ is a smooth Segre surface.
If $e_1,e_2,e_3$ and $e_4$ are the exceptional curves of the blowup,
the following 16 curves are all $(-1)$-curves on $S$:
\begin{gather}
e_1,e_2,e_3,e_4,\label{exc}\\
(1,0) - e_1,\,\, (1,0) - e_2,\,\, (1,0) - e_3,\,\, (1,0) - e_4,
\label{10}\\
(0,1) - e_1,\,\, (0,1) - e_2,\,\, (0,1) - e_3,\,\, (0,1) - e_4,\label{01}\\
(1,1) - e_{234},\,\, (1,1) - e_{134},\,\, (1,1) - e_{124},\,\, (1,1) - e_{123}, 
\label{11}
\end{gather}
where in the last line, $e_{234}$ means $e_2 + e_3 + e_4$ and so on.
All these curves are mapped to lines in $\CP_4$ isomorphically
by the anti-canonical system.
By adding up all these 16 classes, we obtain the class
\begin{align}\label{brnch}
(8,8) - 4e_{1234}.
\end{align}
By Propositions \ref{prop:line1} and \ref{prop:br3}, this is exactly the class of the branch divisor
of the generically finite double covering $\bm D_S\to S$
for the smooth Segre surface $S$.
The class \eqref{brnch} is exactly the class $4K_S\inv$.
Moreover, since each component of the branch divisor is a line,
all their intersections are transverse.
Therefore, any singularity of the branch divisor is an $A_1$-point.
Hence, any singularity of the double covering branched along all lines on $S$ is an $A_1$-point.
If $\pi_1$ denotes the covering map as before, the canonical class of this double covering is given by
$$
\pi_1^*\left(K_S + \frac12\cdot 4K_S\inv\right) \simeq 
\pi_1^* K_S\inv.
$$
Since $K_S\inv$ is ample,
this means that the canonical class of the double covering is 
big. In particular, we obtain

\begin{proposition}\label{prop:gt}
The cuspidal locus in $S^*$ of any smooth Segre surface $S$ is an irreducible surface of 
general type.
\end{proposition}

\subsection{Singularities of the cuspidal locus}
\label{ss:sing}
Finally, we discuss singularities of the cuspidal locus
in the dual variety $S^*$ when the Segre surface $S$ has a line that does not
go through any singularity of $S$
(namely when the divisor $\bm D_S$ is irreducible).
By Proposition \ref{prop:line1}, such a line $l$ is always 
a simple branch divisor of the generically finite double covering 
$\bm D_S\to S$. In particular, the divisor $\bm D_S$ is
smooth at generic points of the ramification curve over $l$.
In this subsection, we first identify the image of 
this ramification curve under the projection $\pi_2$ in 
a concrete form, and next show that
the image $\pi_2(\bm D_S)$, namely the cuspidal locus in $S^*$,
has ordinary cusps along the image of the ramification curve over $l$.
For these purposes, we first prove a few lemmas.

\begin{lemma}\label{lemma:qdr1}
If a line $l$ on a Segre surface $S$ does not pass any singularity of $S$, then $\pi_1\inv(l)\simeq l\times\CP^1$ holomorphically.
\end{lemma}

\proof
Take any $p\in S_{\reg}$ and let $N_p\simeq\CC^2$ denote the fiber
of the normal bundle $N_{S_{\reg}/\CP_4}\to S_{\reg}$ over the point $p$.
The fiber of the projection 
$\pi_1:I(S)\to S$ over $p$ is identified with $\mathbb P(N_p)\simeq\CP_1$.
Therefore we have a natural isomorphism 
$$\pi_1\inv(l)\simeq\mathbb P\big(N_{S_{\reg}/\CP_4}|_l\big).$$
Take a pair of quadrics $Q_1$ and $Q_2$
in $\CP_4$ which satisfy $Q_1\cap Q_2 = S$.
Then we have an isomorphism $N_{S_{\reg}/\CP_4}\simeq N_{S_{\reg}/Q_1}\oplus N_{S_{\reg}/Q_2}$.
Further, we have $N_{S_{\reg}/Q_1}|_l\simeq [Q_2]|_l\simeq\ms O(2)$
and similarly $N_{S_{\reg}/Q_2}|_l\simeq \ms O(2)$.
Hence,  $N_{S_{\reg}/\CP_4} |_l \simeq \ms O(2)^{\oplus 2}$.
These mean the assertion of the lemma.
\proofend

\begin{lemma}
Let $S$ and $l$ be as in the previous lemma,
and $l^*$ the 2-plane in $\CP_4^*$ formed by 
hyperplanes which contain the line $l$.
Then we have $\pi_2\big(\pi_1\inv(l)\big) = l^*$.
\end{lemma}

\proof
Suppose $(p,H)\in \pi_1\inv(l)$. Then $p\in l$ and $T_pS\subset H$.
As $l\subset T_pS$, these mean $l\subset H$. Namely $H\in l^*$.
Hence $\pi_2\big(\pi_1\inv(l)\big) \subset l^*$.
To show the reverse inclusion, it suffices to show that 
the image $\pi_2\big(\pi_1\inv(l)\big)$ is 2-dimensional
because $l^*$ is irreducible.
If $p$ is any point of $l$, then the fiber $\pi_1\inv(p)$ is identified
with a line in $\CP_4^*$ formed by hyperplanes containing $T_pS$.
So, if the image $\pi_2\big(\pi_1\inv(l)\big)$ would be 1-dimensional,
the tangent plane $T_pS$ has to be independent of a choice of 
a point $p\in l$.
If this is actually the case, any $H$ containing the common 2-plane
has to satisfy $H|_S \ge 2l$. 
But by Lemma \ref{lemma:line1}, such a hyperplane does not exist.
Hence the tangent plane $T_pS$ really varies as a point $p$ moves on $l$.
\proofend

\medskip
The next lemma identifies the map $\pi_1\inv(l)\to l^*$ in a concrete form:

\begin{lemma}\label{lemma:invol}
Let $l\subset S$  be as in Lemma \ref{lemma:qdr1},
so that $\pi_1\inv(l)\simeq\qdr$.
Then the restriction of the projection $\pi_2:I(S)\lras S^*$ to the divisor $\pi_1\inv(l)$
can be identified with a quotient map from $\pi_1\inv(l)$ 
under the involution on $\pi_1\inv(l)\simeq\qdr$ given by 
the reflection $(x,y)\mapsto (y,x)$.
\end{lemma}

\proof
By the previous lemma, the image $\pi_2\big(\pi_1\inv(l)\big)$
is a 2-plane.
It is well-known that any degree-two morphism
from $\qdr$ to $\CP_2$ is identified with
the quotient map under the reflection as in the lemma.
Hence it is enough to show that for a generic
point $H\in l^*$, there exist precisely two points
$p,q\in S_{\reg}$ which satisfy $(p,H)\in \pi_1\inv(l)$
and $(q,H)\in \pi_1\inv(l)$.
The condition $(p,H)\in \pi_1\inv(l)$ is equivalent to the conditions 
$p\in l$ and $T_pS\subset H$,
and therefore, the presence of the above two points $(p,H)$ and $(q,H)$
is equivalent to 
the condition that the hyperplane $H$
is tangent to $S$ at $p$ and $q$.
By \cite[Proposition 3.12]{H20}, if $H$ is generic
in the 2-plane $l^*$, then $H|_S$ is of the form 
$l+C$, where $C$ is a smooth rational cubic curve
intersecting $l$ transversely at exactly two points.
Hence $H$ is indeed tangent to $S$ at two points on $l$.
\proofend

\medskip
By Lemma \ref{lemma:invol}, the ramification divisor 
of the projection $\pi_1\inv(l)\to l^*$ can be identified
with the diagonal of $\pi_1\inv(l)\simeq\qdr$,
and the branch divisor is a smooth conic in the 2-plane $l^*$.
In the sequel, we denote these ramification divisor
and branch divisor by $\DDD$ and $\ol\DDD$ respectively.
The divisor $\bm D_S$ is
smooth at generic points of the ramification divisor $\DDD$.
For a generic point $(p,H)$ belonging to $\DDD$, 
the structure of the hyperplane section $H|_S$ is described as follows.

\begin{lemma}\label{lemma:lC}
%Let $l\subset S$  be as in Lemma \ref{lemma:qdr1}.
%Then for a generic hyperplane $H\subset\CP_4$ which belongs
%to the branch divisor $\ol\DDD\subset l^*=\CP_2$,
The section $H|_S$ is of the form $l + C$,
where $C$ is a rational normal cubic curve which intersects $l$
at a unique point, and which is tangent to $l$ at the point
with order two.
\end{lemma}

\proof
By the two projections $\pi_1$ and $\pi_2$,
there are natural isomorphisms $l\simeq\DDD\simeq\ol\DDD$.
As the conic $\ol\DDD$ is the branch divisor of 
the double covering $\pi_1\inv(l)\to l^*$,
a hyperplane $H$ belongs to $\ol\DDD$ iff the residual cubic curve $C = H|_S - l$ satisfies
the property that $C\cap l$ consists of a single point.
Let $p$ be this point.
Of course this is a smooth point of $S$,
and $p\in l$ and $H\in\ol\DDD$ are identified under the above isomorphism $l\simeq\ol\DDD$.
To show the lemma, it is enough to show that if the point $p$ does not belong to
any other line on $S$, then the residual curve $C$ satisfies 
the properties in the proposition.

For intersection numbers, we have
\begin{align}\label{ctan}
(C,l)_S = (H|_S - l, l)_S = (H,l)_{\CP_4} - (l,l)_S 
= 1 - (-1) = 2.
\end{align}
Since $C\cap l=\{p\}$ as above, this means that 
$C$ and $l$ do not intersect transversely at $p$.
If $C$ would be reducible, then from the above assumption on
the point $p$,
$C$ has to have a smooth conic $D$ as a component
and $D$ has to be tangent to $l$ at $p$.
Then the curve $l+D$ lies on a 2-plane, and this contradicts
Proposition \ref{prop:SI2}.
Hence $C$ is an irreducible cubic curve through $p$
which does not intersect transversely at $p$.
Moreover, again by Proposition \ref{prop:SI2}, $C$ has to be non-planer.
Therefore, $C$ is a non-degenerate irreducible cubic curve
in $\CP_3$. Hence it has to be a rational cubic curve,
and by \eqref{ctan}, it has to be tangent to $l$ at $p$ with order two.
\proofend

\medskip
Thus, the singularity of hyperplane sections of $S$ given by generic points
of the conic $\ol\DDD$ is not an ordinary cusp but a tacnode.
So it would be natural to expect that the cuspidal locus
is singular along the conic $\ol\DDD$.
The next result shows that this is actually the case
and identify the singularity along $\ol\DDD$.

\begin{proposition}\label{prop:DDD}
Let $S\subset\CP_4$ be a Segre surface and $l\subset S$
a line satisfying $l\cap \Sing(S)=\emptyset$.
Then the cuspidal locus in $S^*$ has ordinary cusps
along generic points of the conic $\ol\DDD$
in the 2-plane $l^*$.
\end{proposition}

\proof
%This can be shown in a similar way to Theorem \ref{thm:cuspcusp2}
%and Proposition \ref{prop:vA3}.
Let $H\subset\CP_4$ be a generic hyperplane belonging to the conic $\ol\DDD$, so that by the previous proposition, the section $H|_S$
is of the form $l + C$, where $l$ is a line and $C$ is a rational normal cubic that is tangent to $l$ at a point.
Let $p$ be the tangent point.
Similarly to the situation in Proposition \ref{prop:vA3},
the curve $l+C$ is singular only at $p$, and the type
of the singularity is a tacnode.
On this reducible curve, we still have an exact sequence 
\begin{align}\label{eqnn5}
0\lras N'_{l + C} \lras N_{l + C} \lras T^1_{l + C} \lras 0.
\end{align}
By Proposition \ref{prop:SI1} and Lemma \ref{lemma:lC}, we have $l^2 = -1, C^2 = 1$,
and $l.\,C=2$. Hence we obtain 
\begin{align}\label{NlC}
N_{l + C}|_{l}\simeq [l + C]|_{l}\simeq \ms O_l(1)
\qandq
N_{l + C}|_{C}\simeq [l + C]|_{C}\simeq \ms O_C(3).
\end{align}
Let $\ms J$ be the Jacobian ideal sheaf of the curve $l+C$,
and write $\ms J_{l+C}$ for $\ms J|_{l+C}$.
Let $\nu:l\sqcup C\to l+C$ be the normalization of $l+C$.
Then the calculations in the proof of Proposition \ref{prop:vA3}
works without any change, and as in \eqref{nuJ2}, we obtain
$$
\big(\nu^*\ms J_{l+C}\big)|_C\simeq \ms O_C(-2)
\qandq
\big(\nu^*\ms J_{l+C}\big)|_l\simeq \ms O_l(-2).
$$
From these and \eqref{NlC}, we obtain 
$$
\nu^* N'_{l+C} \simeq \ms O_l(-1) \sqcup \ms O_C(1),
$$
and hence, also an exact sequence
\begin{align}\label{ses:Ls2}
0 \lras  N'_{l + C} \lras \nu_*\big(\ms O_l(-1) \sqcup \ms O_C(1)\big)
\lras \CC_p \lras 0.
\end{align}
From this, we again have
$$
H^0\big(N'_{l+C} \big)\simeq\CC \qandq H^1\big(N'_{l + C} \big) = 0.
$$
Hence we again obtain that 
any first order displacement of
the curve $l+C$ in $S$ is unobstructed, and 
that the versal family of 
equi-singular displacements of $C_1+C_2$ in $S$ is 1-dimensional.
From \eqref{eqnn5}, we also obtain an exact sequence
\begin{align}\label{eqnn7}
0\lras H^0\big(N'_{l+C}\big) \lras H^0\big(N_{l+C}\big) \lras 
H^0\big(T^1_{l+C}\big) \lras 0.
\end{align}

For the versal deformation of the tacnode,
as discussed right before the proof of Theorem \ref{thm:cuspcusp2},
the versal family is smooth and 3-dimensional.
In the notation and argument there, 
the fiber of 
the lift of the versal family by the Galois cover $(t_1,t_2,t_3,t_4) \mapsto (s_1,s_2,s_3)$
%with the constraint $t_1 + t_2 + t_3 + t_4=0$
has an $A_2$-singularity iff $(t_1,t_2,t_3,t_4)$ is equal to 
$$(t,t,t,-3t),\,\, (t,t,-3t,t),\,\, (t,-3t,t,t) {\text{ or }} (-3t,t,t,t)$$
for some $t\neq 0$.
Therefore, the fiber of the versal family over a point $(s_1,s_2,s_3)$ has an $A_2$-singularity iff
it is of the form $(-6t^2, -8t^3, -3t^4)$ for some $t\neq 0$.
Similarly to the notations we used in Section \ref{ss:cuspidal2}, 
let $\ms A_2\subset\CC^3$ be the locus formed by points in this form.
Then the closure of $\ms A_2$ has an $A_2$-singularity
at the origin,
since, up to a non-zero constant multiple to each factor, 
the map $(z_1,z_2)\longmapsto (z_1,z_2,z_1^2)$
gives an isomorphism from a curve with an ordinary cusp
to the above closure.

As in the proof of Theorem \ref{thm:cuspcusp2},
take a small neighborhood $B$ of the point $H\in \DDD$ in $\CP_4^*$,
and regard $B$ as a parameter space for displacements of the 
hyperplane section $l+C$ in $S$.
Then by versality, we have an induced holomorphic map $f:B \to \CC^3$,
and the differential $(df)_H$ is identified with 
the map $H^0\big(N_{l+C}\big) \lras 
H^0\big(T^1_{l+C}\big)$ in \eqref{eqnn7}.
The intersection of $B$ with the cuspidal locus in $S^*$
is exactly $f\inv(\ms A_2\cup\{H\})$.
By the surjectivity of the map $H^0\big(N_{l+C}\big) \lras 
H^0\big(T^1_{l+C}\big)$ and the above property about the closure of
the locus $\ms A_2$, 
the cuspidal locus has $A_2$-singularities along 
$f\inv(0)$.
Since sections of $S$ by hyperplanes
belonging to the conic $\ol\DDD$ have tacnode by Lemma \ref{lemma:lC}, 
$f\inv(0)$ is equal to $\ol\DDD\cap B$.
Thus, the cuspidal locus has $A_2$-singularities along generic
points of $\ol\DDD$.
\proofend

\medskip
Proposition \ref{prop:DDD} in particular means that,
for a Segre surface $S$ which has a line not passing any singularity of $S$,
the cuspidal locus in $S^*$ is non-normal, and that
the projection $\pi_2:\bm D_S\to \pi_2(\bm D_S)$ gives the normalization along the conic $\ol\DDD$.

\section{Appendix. explicit calculations for the cuspidal locus}
\label{s:scl2}

Here is a list of Segre symbols for
all Segre surfaces that have a line 
which passes two singularities of the surfaces:
$$
\begin{array}{cccccccc}%\label{}
{\text{symbol}} & [(11)(11)1] & [2(11)1] & [221] & [3(11)]  &  [32]  &  [(21)(11)]  &  [2(21)] \\
{\text{singularities}} & 4A_1 & 3A_1 & 2A_1 & A_2 + 2A_1  &  A_1+A_2  &  A_3+2A_1  &  A_1+A_3 
\end{array}
$$
A Segre surface with symbol $[(11)111]$
has two singularities, but it is not included 
in this list because the line through these points
does not lie on the surface.
%All lines connecting $A_1$-point and non-$A_1$-point
%are contained in the Segre surfaces.

In this subsection, to each surface in the above list,
by making use of the common tangent plane
to the surface along the line $l$ connecting two singularities of $S$
(see Lemma \ref{lemma:line3}), we calculate the two functions
$F$ and $G$ in concrete forms to some extent.
Next, we determine the precise multiplicity of the divisor
$\pi_1\inv(l)$ as a component of $\bm D_1$,
and show that the line $l$ is always {\em not} a branch divisor of 
the generically finite double covering $\bm D_S\to S$.
Applying this to the cases where the Segre symbol is 
 $[32]$ or $[221]$, this in particular proves Proposition \ref{prop:br3} 
about the branch divisor of the covering  $\bm D_S\to S$.

In the following argument, $(X_0,X_1,X_2,X_3,X_4)$
are homogeneous coordinates on $\CP_4$.
For brevity, we use the symbol $e_i$ ($0\le i\le 4)$ to mean
the point of $\CP_4$ whose unique non-zero entry is  the $i$-th one.
For given two functions $F(x,y)$ and $G(x,y)$, we put 
$$
K := F_{xx}G_{yy} + F_{yy}G_{xx} - 2F_{xy}G_{xy}
$$
for simplicity, so that the defining equation of 
the divisor $\bm D=\bm D_1 + \bm D_S$ in the incidence variety $I(S)$
is given by
$$
{\bf H} = \Hess (F) \lmd^2 + K \lmd \mu + \Hess(G) \mu^2.
$$
The way how one can derive the normal forms for the pair of the defining
quadratic polynomials for each type of Segre surface is
briefly explained in \cite[Section 3.2]{H20}.

\subsection{}
First we discuss a Segre surface $S$ whose symbol is $[(11)(11)1]$.
This would be the easiest case in actual calculations.
%This surface is a toric surface, 
Let $\aaa,\bbb,\ccc$ be mutually different complex numbers. The normal forms for the equations of $S$ on $\CP_4$ are given by 
$$
\aaa X_0^2  + \aaa X_1 ^2 + \bbb X_2^2  + \bbb X_3^2 + \ccc X_4^2
= X_0^2 +  X_1^2 + X_2^2 + X_3^2 + X_4^2 = 0.
$$
It is more convenient to make an obvious coordinate change
and rewrite the equations as
%$$
%Y_0 = X_0 + iX_1, \,\,Y_1 = X_0 - iX_1, \,\,
%Y_2 = X_2 + iX_3, \,\,Y_3 = X_2 - iX_3, \,\,
%Y_4 = X_4.
%$$
%Then the two equations become
$$
\aaa X_0 X_1 + \bbb X_2 X_3 + \ccc X_4^2 = 
X_0 X_1 +  X_2 X_3 +  X_4^2 = 0.
$$
The surface $S$ has exactly four singularities, and all of them 
are $A_1$-points.
They are concretely given by 
$
e_i
$ with  $i = 0,1,2,3$.
The involution of $\CP_4$ which exchanges $X_0$ and $X_1$
keeps $S$ invariant, and it exchanges
the two $A_1$-points $e_0$ and $e_1$.
Similarly, the involution of $\CP_4$ which exchanges $X_2$ and $X_3$
keeps $S$ invariant, and it exchanges
the two points $e_2$ and $e_3$.
Among six lines connecting the four $A_1$-points,  
the four lines
$\ol{e_0e_2}, \ol{e_0e_3}, \ol{e_1e_2}$ and $\ol{e_1e_3}$ lie
on $S$.
%\begin{align*}
%\ol{e_0 e_2} = \{ X_1  = X_2 - i X_3 = X_4 = 0\},\,\,
%\ol{e_2e_4} = \{ X_0 + i X_1 = X_2 + i X_3 = X_4 = 0\},\\
%\ol{e_1e_4} = \{ X_0 - i X_1 = X_2 + i X_3 = X_4 = 0\},\,\,
%\ol{e_2e_3} = \{ X_0 + i X_1 = X_2 - i X_3 = X_4 = 0\}.
%\end{align*}
The group generated by the above two involutions on $S$
acts transitively on the set of these four lines.
So in the following, we choose a line $l:=\ol{e_0e_2}=\{X_1 = X_3 = X_4=0\}$
and calculate an equation of the divisor $\bm D$ in the incidence variety $I(S)$,
in a neighborhood of points of $\pi_1\inv(l)$.
The common tangent plane to $S$ along points on $l\minus\{e_0,e_2\}$ is given by 
\begin{align}\label{T61}
T = \{X_1 = X_3 =0\}.
\end{align}

%Similarly, for other three lines, we have
%$$
%T_{12}= \{ X_0 = X_3 =0\},\quad
%T_{03}= \{ X_1 = X_2 =0\},\quad
%T_{13}= \{ X_0 = X_2 =0\},
%$$
In the sequel, we work on the open subset $\{X_0\neq 0\}=\CC^4$ and use
$x _ i = X_i/X_0$, $1\le i\le 4$, as coordinates on it.
Then from \eqref{T61}, we may use
$(x,y) := (x_2,x_4)$ as coordinates on $S$ around any point
of $l\minus\{e_0,e_2\}$.
The two points $e_0$ and $e_2$ are excluded
since these are singular points of $S$.
On $S$, we have $l = \{y=0\}$.
From the second equation of $S$, we obtain
\begin{align}\label{G7}
x_1 = -x x_3 - y^2. 
\end{align}
Substituting this into the affine form of the first equation of $S$, we obtain
\begin{align}\label{quadF6}
(\bbb-\aaa) x x_3 + (\ccc - \aaa) y^2 = 0.
\end{align}
%If we put $ x = 0$ in this equation, 
%we obtain $(\ccc-\aaa)y^2 =0$.
%Since $\aaa\neq \ccc$, this mean that the hyperplane 
%section by $ X_2 = 0$ contains the line $l$
%by multiplicity precisely two.
We put
$$
F = x_3 \qandq
G = x_1.
$$
From \eqref{quadF6} and \eqref{G7}, we obtain
$$
F =  \frac{\aaa - \ccc}{\bbb - \aaa}\,\frac{y^2}x,\quad
G = \frac{\ccc - \bbb}{\bbb - \aaa} y^2.
$$
Thus we have obtained the two functions $F$ and $G$ 
in explicit forms.
(The divisibility of $F$ and $G$ by $y^2$ follows from 
Lemma \ref{lemma:line3} from the beginning.)
From these, we can compute 
$$
\Hess(F) = \Hess (G) = 0,\quad
K = F_{xx}G_{yy} = \frac{4 (\aaa - \ccc) (\ccc -\bbb)}
{(\bbb - \aaa)^2}\,\frac{y^2}{x^3}.
$$
Hence, we obtain 
$$
{\bf H} = K\lmd\mu=
\frac{4 (\aaa - \ccc) (\ccc -\bbb)}
{(\bbb - \aaa)^2}\,\frac{y^2}{x^3}\lmd\mu.
$$
Since this is divisible precisely by $y^2$,
we obtain that 
{\em the component $\pi_1\inv(l)$ is included
in $\bm D$ and $\bm D_1$ with multiplicity precisely two.}
Also we obtain that the equation of the divisor $\bm D_S = \bm D - \bm D_1$
is simply $\lmd\mu=0$. 
This implies that {\em the generically finite double covering $\bm D_S\to S$ 
does not have the line $l$ as a branch divisor.}

\subsection{}
Next, we discuss a Segre surface $S$ whose symbol is $[2(11)1]$.   
After making a simple coordinate change as in the previous case,
equations of $S$ are given by 
$$
2\aaa X_0 X_1 + X_1 ^2 + \bbb X_2 X_3 + \ccc X_4^2
=2 X_0 X_1 + X_2 X_3 + X_4^2 = 0,
$$
where $\aaa,\bbb,\ccc$ are again mutually different complex numbers.
%$$
%2\aaa X_0 X_1 + X_1 ^2 + \bbb X_2^2  + \bbb X_3^2 + \ccc X_4^2
%=2 X_0 X_1 + X_2^2 + X_3^2 + X_4^2 = 0.
%$$
This surface has exactly three singularities, 
and they are the points $e_0,e_2$ and $e_3$.
All of them are $A_1$-points.
%They are concretely given by 
%$$
%e_0,\,\, p_1:= (0,0,1,-i,0) \qandq p_1:= (0,0,1,i,0).
%$$
The involution of $\CP_4$ which exchanges $X_2$ and $X_3$
keeps $S$ invariant, and it exchanges $e_2$ and $e_3$.
The two lines $\ol{e_0e_2}%= \{X_1 = X_3 = X_4 = 0\}
$ and $\ol{e_0e_3}%= \{X_1 = X_2 = X_4 = 0\}
$ are contained in $S$, and the line $\ol{e_2e_3}$ is not.
The first two lines are exchanged by the above involution.
In the following, we choose the line $l:= \ol{e_0 e_2}= \{X_1 = X_3 = X_4 = 0\}$ and calculate
an equation of the divisor $\bm D$ around points over 
$l\minus\{e_0,e_2\}$.
The common tangent plane to $S$ along points on $l\minus\{e_0,e_2\}$ is given by 
\begin{align}\label{T62}
T = \{ X_1 =  X_3 = 0\}.
\end{align}

In the sequel, we again work on the open subset $\{X_0\neq 0\}=\CC^4$ and use
$x _ i = X_i/X_0$, $1\le i\le 4$, as coordinates on it.
We put $(x,y):=(x_2,x_4)$ and use these as local coordinates
on $S$ around any point on $l\minus\{e_0,e_2\}$.
On $S$, the line $l$ is again defined by $y=0$. 
From the affine form of the second equation, we obtain
\begin{align}\label{G4}
x_1 = -\frac12 \left(
 x x_3 + y^2\right).
\end{align}
Substituting this into the affine form of the first equation, we obtain
\begin{align}\label{quadF3}
x ^2 x_3 ^2
+ 2 x\left\{ 2(\bbb - \aaa) + y^2 \right\}  x_3   +  y ^2 \left\{ y ^2 + 4(\ccc-\aaa) \right\} = 0.
\end{align}
%If we put $ x = 0$ in this equation, 
%we obtain $y^2 \left\{ y^2 - 4(\ccc-\aaa) \right\}=0$.
%Since $\aaa\neq \ccc$, this mean that the hyperplane 
%section by $ X_2  = 0$ contains the line $l$
%by multiplicity precisely two.
We put
$$
F = x_3, \quad
G = -2 x_1.
$$
Then \eqref{quadF3} and \eqref{G4} can be written respectively as 
\begin{align}\label{quadF4}
 x^2 F^2 + 2 x\left\{ 2(\bbb - \aaa) + y^2\right\}  F
+ y ^2 \left\{ y ^2 + 4(\ccc-\aaa) \right\} = 0 
\end{align}
and
\begin{align}\label{G5}
G =  x F + y^2 .
\end{align}
Any hyperplane section
of $S$ containing the common tangent plane 
\eqref{T62} includes $2l$ as a subdivisor,
so we may put $F = y^2 f$ and $G= y^2g$
for some holomorphic functions $f$ and $g$. 
Substituting these
into \eqref{quadF4} and \eqref{G5} and dividing by $y^2$, 
we obtain 
\begin{align}\label{quadF5}
 x^2 y^2f ^2 + 2 x\left\{ 2(\bbb - \aaa) + y^2\right\} f
+ \left\{ y ^2 + 4(\ccc-\aaa) \right\} = 0 
\end{align}
and
\begin{align}\label{G6}
g = x f + 1.
\end{align}
Differentiating \eqref{quadF5} by $x$ and using $\aaa\neq\bbb$, we readily obtain 
a divisibility
\begin{align}\label{div1}
y \set (f + xf_x).
\end{align}
If we write $f$ as 
$f = u(x) + y v(x,y)$, we have
$$
f + xf_x = (u +xu') + y (v + xv_x).
$$
From \eqref{div1}, this in particular means $u + xu'=0$.
Since the function $f(x,y)$ is known to be holomorphic
only in a neighborhood of a point on $l$ which is different
from singularities of $S$, and since the singularity is 
the origin $(x,y) = (0,0)$ in the present coordinates,
we have to allow a pole for $u=u(x)$ at $x=0$.
From the last equation $u + xu'=0$,
this means $u(x) = \frac cx$ for some constant $c$.
So we have
\begin{align}\label{f}
f = \frac cx + yv(x,y).
\end{align}
If $c=0$, we would obtain $y\set f$, which means $y^3 \set F$.
This implies that the section of $S$
by a hyperplane $H=\{ X_3=0\}$ has a triple line $3l$ as component.
But this cannot happen because both of the singularities
on $l$ are $A_1$-points, and from Lemma \ref{lemma:line4},
no hyperplane section contains the triple line $3l$.
Hence, $c\neq 0$ holds.
Differentiating \eqref{f} twice by $x$, we obtain 
\begin{align}\label{fxx5}
f_{xx} = \frac{2c}{x^3} + yv_{xx}. %\quad c\neq 0.
\end{align}
This means that $f_{xx}$ is not divisible by $y$. We will soon use these
to determine exact multiplicity of the component 
$\pi_1\inv(l)$ in the divisor $\bm D_1$.

From $F = y^2f$, we easily obtain
\begin{align}\label{HessF5}
\Hess(F) &= y^2 \big\{
f_{xx} (2f + 4yf_y + y^2 f_{yy})
- (2f_x+yf_{xy})^2
\big\},
%\Hess(G) &= y^2 \big\{
%g_{xx} (2g + 4yg_y + y^2 g_{yy})
%- (2g_x+yg_{xy})^2\label{HessG1}
%\big\},
\end{align}
and similarly for $\Hess(G)$,
meaning $y^2 \set \Hess(F)$ and $y^2\set \Hess(G)$.
From \eqref{HessF5}, it follows that 
$$y^3\set \Hess(F) \Longleftrightarrow
y\set (ff_{xx} -2f_x^2).$$
%and similarly for $\Hess(G)$.
From \eqref{f} and \eqref{fxx5}, we readily obtain 
$$ff_{xx}\equiv \frac{2c^2}{x^4} \mod y
\qandq
f_x^2 \equiv \frac{c^2}{x^4} \mod y.
$$
These mean
$y\set (ff_{xx} -2f_x^2)$.
Hence we obtain $y^3\set \Hess(F)$.

Recall that we have written $F_{xx}G_{yy} + F_{yy}G_{xx} - 2F_{xy}G_{xy}$ as $K$.
Next we show $y^2 \set K$ and $ y^3 \not| K.$
We have $F_{xx} = y^2f_{xx}$, and 
since $y\not|\,f_{xx}$ from \eqref{fxx5}, we have 
$y^2\set F_{xx}$ and $y^3\not|\,F_{xx}$.
From \eqref{G6} and \eqref{f}, we readily obtain 
\begin{align}\label{g50}
g = c + xyv + 1.
\end{align}
If $c=-1$, then $y\set g$ follows
and hence $y^3\set G$. 
Again this contradicts Lemma \ref{lemma:line4}.
Hence $c\neq -1$.
The constant term of $G_{yy}$ is easily seen to be $2(c+1)\neq 0$.
From these, we obtain $y^2\set F_{xx}G_{yy}$ and $y^3\not|\, F_{xx}G_{yy}$.
Also, from \eqref{g50}, we easily obtain $y\set g_{xx}$.
This means $y^3\set G_{xx}$ and hence $y^3\set F_{yy}G_{xx}$.
Further, it is immediate to see $y\set F_{xy}$,
and from \eqref{g50} we also have $y^2 \set G_{xy}$.
Hence $y^3 \set F_{xy}G_{xy}$. So $y\not\set G_{yy}$.
From these, we obtain $y^2 \set K$ and $ y^3 \not| K.$

Thus, we obtain $y^2 \set \bf H$ and $y^3\not|\, \bf H$.
Hence, {\em the multiplicity of the divisor $\pi_1\inv(l)$ in 
$\bm D_1$ is precisely two}.
The discriminant of the quadratic polynomial ${\bf H}/y^2$ is equal to
\begin{align}\label{Hess13}
\left(\frac{K}{y^2}\right)^2 - 
4\frac{\Hess(F)}{y^2}\frac{\Hess(G)}{y^2}. 
\end{align}
The second term of this is divisible by $y$ since $y^3\set \Hess(F)$ as above. But the first term is not divisible by $y$ since
$y^3\not| K$ as above.
Hence \eqref{Hess13} is not divisible by $y$.
Thus we can conclude that {\em the line $l$ is not a component of 
the branch locus of 
the generically finite double covering
$\bm D_S \to S$.}

%This means 
%$y\set g_x$ and $y\set g_{xx}$.
%Hence
%$
%y\set (gg_{xx} -2g_x^2).
%$
%Therefore $y^3\set \Hess(G)$.
%Also, we obtain from \eqref{f} and \eqref{G6}  that 
%$$y^2 \set K \qandq  y^3 \not| K.$$
%

\subsection{}
Next we discuss the case where the Segre symbol of $S$ is $[221]$. The normal forms of the equations of $S$ are given by 
$$
2\aaa X_0 X_1 + X_1 ^2 + 2 \bbb X_2 X_3 + X_3^2 + \ccc X_4^2
=2 X_0 X_1 + 2 X_2 X_3 + X_4^2=0,
$$
where $\aaa,\bbb,\ccc$ are mutually different complex numbers. 
The two points $e_0$ and $e_2$ are $A_1$-singularities of $S$,
and these are all singularities of $S$.
The line connecting these two singularities is
$$
l = \{X_1 = X_3 = X_4 = 0\},
$$
and this is contained in $S$. 
The common tangent plane to $S$ along points on $l\minus\{e_0,e_2\}$ is again given by 
$$
T = \{X_1 = X_3 = 0\}.
$$

In the sequel, we again work on the open subset $\{X_0\neq 0\}=\CC^4$ and use
$x _ i = X_i/X_0$, $1\le i\le 4$, as coordinates on it.
From the equations of $T$, we may use $(x,y):=(x_2,x_4)$ as local coordinates
on $S$ around any point of $l\minus\{e_0,e_2\}.$
On $S$, the line $l$ is defined by  $y=0$.
From the affine form of the second equation of $S$, we obtain
\begin{align}\label{G3}
x_ 1= -x x_3 - \frac12 y ^2.
\end{align}
Substituting this into the affine form of the first equation,
we readily obtain 
\begin{align}\label{quadF1}
x^2 x_3^2 + x\{y^2 - 2(\aaa-\bbb)\} x_3 
+ \frac 14 y^4 + (\ccc-\aaa) y^2 = 0.
\end{align}
%If we put $x_3=0$ in this equation, we obtain 
%$y^2 ( y^2 + 4\ccc - 4\aaa) = 0$.
%Noting $\aaa\neq \ccc$, this means that 
%the hyperplane section $\{X_3=0\}|_S$ contains the double line $2l$.
We put 
$$
F = x_3 \qandq G = x_1
$$
for the normal directions to $S$.
As in the last case of $[2(11)1]$, 
we may put $F = y^2 f$ and $G= y^2g$.
Then \eqref{quadF1} is divisible by $y^2$, and we obtain that 
the function $f$ is subject to the equation
\begin{align}\label{quadF2}
x^2 y^2 f^2  + x\{y^2 - 2(\aaa-\bbb)\} f
+ \frac 14 y^2 + (\ccc-\aaa)  = 0.
\end{align}
Differentiating this by $x$ and using $\aaa\neq\bbb$, we readily obtain the divisibility
$$
y^2 \set (f + xf_x).
$$
In the same way to the last case of $[2(11)1]$, this means 
\begin{align}\label{f1}
f = \frac cx + yv(x,y),
\end{align}
for some constant $c$ and a holomorphic function $v$.
If $c=0$, we obtain $y^3\set F$, which contradicts 
Lemma \ref{lemma:line4}. So $c\neq 0$.
In the same way to the last case of $[2(11)1]$, we obtain from $F= y^2 f$ that $y^2\set \Hess(F)$, and further, 
%\begin{align}\label{}
%\Hess(F) &= y^2 \big\{
%f_{xx} (2f + 4yf_y + y^2 f_{yy})
%- (2f_x+yf_{xy})^2\label{HessF}
%\big\},
%%\Hess(G) &= y^2 \big\{
%%g_{xx} (2g + 4yg_y + y^2 g_{yy})
%%- (2g_x+yg_{xy})^2\label{HessG1}
%%\big\},
%\end{align}
%and similarly for $\Hess(G)$.
%From this it follows that 
$$y^3\set \Hess(F) \Longleftrightarrow
y\set (ff_{xx} -2f_x^2).$$
Again from \eqref{f1}, we obtain that 
$f$ satisfies the latter condition.
Therefore $y^3\set \Hess(F)$ holds.
On the other hand, 
from \eqref{G3}, we have
$$
G = - x F - \frac12 y^2,
$$
which means 
\begin{align}\label{g1}
g = -\Big(c+\frac 12\Big) - xy v.
\end{align}
This implies $y\set g_x$ and $y\set g_{xx}$.
Hence $y\set (gg_{xx} -2g_x^2)$ also holds.
So, in the same way for $F$, we have $y^3 \set \Hess(G)$.

From $F = y^2f$ we have $F_{xx} = y^2 f_{xx}$,
and $f_{xx}\equiv \frac {2c}{x^3} \mod y$ from \eqref{f1}.
Also, from $G = y^2g$, we have $G_{yy} = 
2g + 4y g_y + y^2 g_{yy}$.
From these, we have 
$$
F_{xx}G_{yy} \equiv \frac{4cg}{x^3}y^2  \mod y^3.
$$
In particular, $y^2\set F_{xx}G_{yy}$ and 
$y^3\not|\, F_{xx}G_{yy}$ as $c\neq 0$.
Next, from $G = y^2g$, we have $G_{xx} = y^2g_{xx}$,
and since $y\set g_{xx}$ as above, we have 
$y^3\set G_{xx}$.
Therefore $y^3 \set F_{yy}G_{xx}$.
Next, $F_{xy} = y(2f_x + y f_{xy})$ and therefore 
$y\set F_{xy}$.
Further, $G=y^2g$ means $G_x = y^2g_x$ and $y\set g_x$
as above. So $y^3\set G_x$,
which means $y^2 \set G_{xy}$.
Therefore $y^3 \set F_{xy}G_{xy}$.
From these, we conclude 
%$$
%K = F_{xx}G_{yy} + F_{yy}G_{xx} - 2F_{xy} \equiv \frac {4cg}{x^3}y^2 \mod y^3.
%$$
%In particular, 
$y^2\set K$ and $y^3\not|\,K$.

%Further we have $y\not| g$ since if so $y^3\set G$ and 
%this contradicts Lemma \ref{lemma:line4}.
These imply that the defining function $\bf H$
of the divisor $\bm D$ satisfies $y^2\set \bf H$ 
and $y^3 \not|\, \bf H$.
Hence, {\em the multiplicity of the divisor $\pi_1\inv(l)$ in 
$\bm D_1$ is precisely two}.
The discriminant of the quadratic polynomial ${\bf H}/y^2$ is 
\begin{align}\label{Hess12}
\left(\frac{K}{y^2}\right)^2 - 
4\frac{\Hess(F)}{y^2}\frac{\Hess(G)}{y^2}.
\end{align}
The second term is divisible by $y$ as $y^3 \set \Hess(F)$.
But the first term is not divisible by $y$
as $K$ is not divisible by $y^3$ as above.
Hence \eqref{Hess12} is not divisible by $y$.
Thus we can conclude that {\em the line $l$ is not a 
branch divisor of 
the generically finite double covering
$\bm D_S \to S$.}

\subsection{}
Next, we discuss a Segre surface $S$ whose symbol is $[3(11)]$.
%This is of type $2A_1 + A_3$.
After making a simple coordinate change, equations of $S$ are given by 
$$
2\aaa X_0 X_2 + \aaa X_1^2 + 2X_1X_2 + \bbb X_3 X_4
= 2X_0 X_2 +  X_1^2 +  X_3 X_4 =0,
%
%2\aaa X_0 X_2 + \aaa X_1^2 + 2X_1X_2 + \bbb X_3^2 + \bbb X_4^2
%= 2X_0 X_2 +  X_1^2 +  X_3^2 +  X_4^2=0,
$$
where $\aaa\neq\bbb$.
This surface has exactly three singularities, and 
they are
$e_0, e_3$ and $e_4$.
The point $e_0$ is an $A_2$-point of $S$, and both $e_3$ and $e_4$ are 
$A_1$-points of $S$.
The two lines $\ol{e_0e_3}$ and $\ol{e_0e_4}$ are contained in $S$,
and $\ol{e_3e_4}$ is not.
The involution of $\CP_4$ which exchanges $X_3$ and $X_4$
keeps $S$ invariant.
It fixes the point $e_0$, and exchanges $e_3$ and $e_4$.
Hence it exchanges the two lines $\ol{e_0e_3}$ and $\ol{e_0e_4}$.
So in the following, we only consider the line 
$$l:= \ol{e_0e_3} =\{X_1 = X_2 = X_4 = 0\}.$$ 
The common tangent plane 
to $S$ at smooth points on this line is given by 
\begin{align}\label{ctp}
T= \{ X_ 2 = X_4 = 0\}.
\end{align}

In the sequel, we again work on the affine open subset
$\{X_0\neq 0\}=\CC^4$ in $\CP_4$, and use $x_i = X_i / X_0$, $1\le i\le 4$, as
coordinates on the open subset. 
From \eqref{ctp},
we may use $(x,y):= (x_3,x_1)$ as coordinates on $S$
around any point of $l\minus\{e_0,e_3\}$.
On $S$, the line $l$ is defined by $y=0$.
The affine form of the second 
defining equation of $S$ is given by
\begin{align}\label{G1}
x_2 = -\frac12\left( y ^2 + x x_4\right).
\end{align}
Substituting this into the affine form of the first defining equation of $S$,
we obtain 
\begin{align}\label{cubic1}
y ^ 3 + (y+\aaa - \bbb) x x_4 = 0.
\end{align}
Thus the surface $S$ is locally identified with the surface in $\CC^3$
defined by this cubic equation.
If we substitute $x_4 = 0$ to this equation,
we obtain $y^3 = 0$.
This means that if $H$ is a hyperplane defined by
$X_4 =0$, then the section $H|_S$ contains 
the triple line $3l$.
This is the unique hyperplane given in Lemma
\ref{lemma:line4} for the present surface $S$.
We put
$$
F = x_4 \qandq
G = -2 x_2,
$$
and we also write $\ddd = \aaa - \bbb$ for short.
Note that $\ddd\neq 0$ since $\aaa\neq\bbb$.
Then from \eqref{cubic1}, the function $F$ is subject to the equation
$$
y^3 + (y + \ddd) x F = 0
\quad{\text{i.e.}}\quad 
F = - \frac{y^3}{x(y+\ddd)}.
$$
Hence from \eqref{G1}, we also obtain
\begin{align}\label{G18}
G = \frac{\ddd y^2}{y+\ddd}.
\end{align}
Thus we were able to obtain the functions $F$ and $G$
in explicit forms.
From this, we obtain 
$$
F_{xx} = -\frac {2y^3} {x^3(y+\ddd)},\quad
y^2 \set F_{xy},\quad %= \frac {y^2} {x^2(y+\ddd)^2},\quad
y \set F_{yy}.% = \frac {y \{4y^3 + 6(\ddd-1) y^2 + (2\ddd-7)y + \ddd^2\}} {x(y+\ddd)^3}.
$$
These imply that $\Hess(F)$ is divisible by $y^4$. (We just need $y^3\set F$ for obtaining this.)
On the other hand, we obtain from \eqref{G18} that 
$$
G_{xx} = G_{xy} = 0
\qandq G_{yy} = \frac{2\ddd^3}{(y+\ddd)^3}.
$$
Hence $\Hess (G) = 0$.
Further we obtain
$$
K=F_{xx}G_{yy} + F_{yy}G_{xx} - 2F_{xy}G_{xy} =
- \frac {4\ddd^3 y^3}{x^3(y+\ddd)^4}.
$$
Using $\ddd\neq 0$, these mean $y^3\set K$ and
$y^4\not| K$.
Therefore, we obtain
$$
{\bf H} =
y^3\lmd\left(\frac{\Hess(F)}{y^3} \lmd - 
\frac {4\ddd^3 }{x^3(y+\ddd)^4}\mu\right).
$$
Noting that $y^4\set \Hess(F)$ and $\ddd\neq 0$ as above,
it follows that {\em the divisor $\pi_1\inv(l)$ is included
as a component of $\bm D_1$ 
with multiplicity precisely three,
and the line $l$ is not contained in the branch divisor
of the generically finite double covering
$\bm D_S \to S$.}

\subsection{}
Next we discuss the case where the symbol of a 
Segre surface $S$ is $[32]$.
%We have $A_1 + A_2$.
The normal forms of defining equations of $S$ are given by 
$$
2\aaa X_0 X_2 + \aaa X_1^2 + 2 X_1 X_2 + 2\bbb X_3X_4 + X_4^2 =  
2 X_0X_2 + X_1^2 + 2 X_3 X_4 = 0.
$$
Singularities of $S$ are the point $e_0$, which is an $A_2$-point, and $e_3$, which is an $A_1$-point.
The line 
$$
l := \{ X_1 = X_2 = X_4 = 0\}
$$
passes $e_0$ and $e_3$ and is contained in $S$.
The common tangent plane to $S$ along points of $l\minus\{e_0,e_3\}$ is given by 
\begin{align}\label{ctp2}
T=\{ X_2 = X_4 = 0\}.
\end{align}

Again we work on the affine open subset $\{X_0\neq 0\}$
and use $x_i=X_i/X_0$ as coordinates on it.
The singular point $e_0$ is the origin.
From \eqref{ctp2}, we may use $(x,y):=(x_3,x_1)$ as coordinates
on $S$ around points on $l\minus\{e_0,e_3\}$.
On $S$ we have $l = \{y=0\}$.
The second equation for $S$ becomes
\begin{align}\label{G2}
2x_2 = - y^2 - 2 x x_4.
\end{align}
Substituting this into the affine form of the first equation of $S$, we obtain an equation
\begin{align}\label{quad1}
x_4^2 - 2x( y +\ddd)  x_4 - y^3 = 0,
\end{align}
where we again put $\ddd=\aaa-\bbb\neq0$.
The surface $S$ is locally identified with this surface in $\CC^3$. 
If we let $x_4 = 0$ in this equation,
we obtain $y^3 = 0$.
This shows that if $H$ is a hyperplane defined by
$X_4 =0$, then the section $H|_S$ contains 
the triple line $3l$.
This is the unique hyperplane given in Lemma
\ref{lemma:line4} in the present surface.
By putting 
$$
F=x_4,\quad G = -2 x_2,
$$
the equations \eqref{quad1} and \eqref{G2}
can be written respectively as
\begin{align}\label{quad2G3}
F ^2 - 2 x (y + \ddd) F - y^3 = 0
\qandq G =  y^2 + 2 x F.
\end{align}
These represent the functions $F$ and $G$ in implicit forms.
Around any point on the line $l$ except the two singularities
of $S$, we may put $F = y^3 f$ and $G = y^2 g$
for some holomorphic functions $f$ and $g$.
Then from \eqref{quad1} and \eqref{G2} we obtain 
\begin{align}\label{quad3G4}
y^3 f^2 - 2 x (y + \ddd) f - 1 = 0
\qandq
g = 1 + 2xy f.
\end{align}
Differentiating the former equation by $x$, we obtain 
\begin{align}\label{quad4}
y^3ff_x - 2(y+\ddd)(f+xf_x) = 0.
\end{align}
Since $\ddd\neq 0$, this means $y^3 \set (f + x f_x)$.
Then in the same way to \eqref{div1}--\eqref{fxx5},
by putting $f = u(x) + y v(x,y)$,
we obtain 
\begin{align}\label{ffxx}
f = \frac cx + yv,\quad
f_{xx} = \frac{2c}{x^3} + yv_{xx}.
\end{align}
If $c=0$, we would obtain $y\set f$, which means $y^4 \set F$.
This implies $\{X_4=0\}|_S = 4l$, which contradicts Lemma \ref{lemma:line4}.
Hence $c\neq 0$.
Then the latter of \eqref{ffxx} means that $f_{xx}$ is not divisible by $y$. 

From $F = y^3 f$, we obtain 
\begin{align}\label{Fx1}
y^3\set F_{xx}, \,\,y^2\set F_{xy}\qandq y\set F_{yy}.
\end{align}
So $y^4 \set \Hess(F)$.
Also, in the same way to \eqref{HessF5},
from $G = y^2g$, we obtain 
\begin{align}\label{}
\Hess(G) = 
%y^2 \big\{
%g_{xx} (2g + 4yg_y + y^2 g_{yy})
%- (2g_x+yg_{xy})^2\notag%
%\big\}\\
y^2 \big\{2(gg_{xx}- 2g_x^2) + yh(x,y)\big\}\label{HessG1}
\end{align}
for some function $h(x,y)$ whose explicit form is not needed.
From the second equation of \eqref{quad3G4}, we obtain 
$y\set g_x$ and $y\set g_{xx}$.
So we have
$$
y \set (gg_{xx}- 2g_x^2).
$$
Therefore from \eqref{HessG1} we obtain $y^3 \set \Hess(G)$.
%We again write $K$ for 
%$F_{xx}G_{yy} + F_{yy}G_{xx} - 2F_{xy}G_{xy}$ for brevity.
We have $F_{xx} = y^3f_{xx}$ from $F= y^3f$,
and $G_{yy}\equiv 2g \mod y$  from $G = y^2 g$.
Hence $F_{xx}G_{yy}\equiv 2 f_{xx} g y^3 \mod y^4$.
As $y \set F_{yy}$ from \eqref{Fx1} and $G_{xx} = y^2g_{xx}$ which is divisible
by $y^3$ as $y\set g_{xx}$, we obtain $y^4 \set  F_{yy}G_{xx}$.
Also, we have $y^2\set F_{xy}$ from \eqref{Fx1}, and $y^2\set G_{xy}$ from
$G_x = y^2g_x, G_{xy} = y (2g_x + yg_{xy})$ and $y\set g_x$.
Hence $y^4 \set F_{xy}G_{xy}$.
From these, we obtain 
$$
K\equiv 2f_{xx}gy^3 \mod y^4.
$$
We have $y \not| f_{xx}$ from the latter of \eqref{ffxx} and $c\neq 0$.
If $g$ would be divisible by $y$, we have $y^3 \set G$, 
which means that the hyperplane section of $S$ by 
$\{X_2 = 0\}$ also contains the triple line $3l$.
This means that for any hyperplane $H$ which contains
the common tangent 2-plane $T$, the section $H|_S$
contains $3l$.
This  cannot happen from Lemma \ref{lemma:line4}.
So $y\not|\, g$.
It follows from these that  $y^4\not| K$.
Therefore, we obtain
$$
y^3\set {\bf H} \qandq y^4 \not|\,{\bf H}.
$$
Hence, {\em the divisor $\pi_1\inv(l)$ is included
as a component of $\bm D_1$ 
with multiplicity precisely three.}

The discriminant of the quadratic polynomial ${\bf H}/y^3$
is equal to
\begin{align}\label{Hess11}
\left(\frac{K}{y^3}\right)^2
-
4\frac{\Hess(F)}{y^3}\frac{\Hess(G)}{y^3} 
\end{align}
Now we have $y\,\big|\, \frac{\Hess(F)}{y^3}$ as
$y^4 \set \Hess(F)$.
So the second term is divisible by $y$.
But the first term is not divisible by $y$
as $K$ is not divisible by $y^4$.
Hence \eqref{Hess11} is not divisible by $y$.
Thus we can again conclude that {\em the line $l$ is not contained 
in the branch locus of 
the generically finite double covering
$\bm D_S \to S$.}

\subsection{}
Next we consider the case where the symbol of a Segre surface $S$
is $[(21)(11)]$.
After making a simple coordinate change to the normal forms,
the defining equations of $S$ are given by 
$$
2\aaa X_0 X_1 +  X_1^2 + \aaa X_2^2 + \bbb X_3X_4 = 
2 X_0 X_1 + X_2^2 +  X_3 X_4 = 0,
$$
where $\aaa\neq \bbb$.
Singularities of $S$ are the point $e_0$, which is an $A_3$-point, and $e_3$ and $e_4$, which are $A_1$-points.
The two lines $\ol{e_0e_3}$ and $\ol{e_0e_4}$ are contained in 
$S$, and the line $\ol{e_3e_4}$ is not.
The involution on $\CP_4$ which exchanges $X_3$ and $X_4$ preserves $S$ and it fixes the $A_3$-point $e_0$.
The $A_1$-points $e_3$ and $e_4$ are exchanged by 
this involution.
Therefore, the two lines $\ol{e_0e_3}$ and $\ol{e_0e_4}$ are
exchanged by the involution.
So in the following, we choose the line 
$$
l := \ol{e_0e_3}=\{ X_1 = X_2 = X_4 = 0\},
$$
and work around points on $l\minus\{e_0,e_3\}$.
The common tangent plane to $S$ along these points is given by 
\begin{align}\label{ctp10}
T=\{ X_1 = X_4 = 0\}.
\end{align}

Again we work on the affine open subset $\{X_0\neq 0\}$
and use $x_i=X_i/X_0$ as coordinates on it.
The singular point $e_0$ is the origin.
From \eqref{ctp10}, we may use $(x,y):=(x_3,x_2)$ as coordinates
on $S$ around points on $l\minus\{e_0,e_3\}$.
We then have $l = \{y=0\}$ on $S$.
The second equation for $S$ becomes
\begin{align}\label{G10}
2 x_1 = - y^2 - x_4 x.
\end{align}
Substituting this into an affine form of the first defining equation of $S$ and again putting $\ddd:=\aaa-\bbb\neq0$, we obtain 
\begin{align}\label{quad10}
(y^2 + xx_4)^2 - 4\ddd xx_4 = 0.
\end{align}
The surface $S$ is locally identified with this surface in $\CC^3$. 
If we put $x_4 = 0$ in this equation, 
one obtains $y^4 =0$.
This shows that 
%if $H$ is a hyperplane defined by
%$X_4 =0$, then 
the hyperplane section $\{X_4=0\}|_S$ contains $4l$.
This is the unique hyperplane given in Lemma 
\ref{lemma:line4} in the present surface,
and we have the coincidence $\{X_4=0\}|_S = 4l$.
By putting 
$$
F=x_4,\quad G = -2x_1,
$$
the equations \eqref{quad10} and \eqref{G10}
can be written respectively as
\begin{align}\label{quad2G30}
x^2 F ^2 + 2 x (y^2 -2 \ddd) F + y^4 = 0
\qandq G =  y^2  + x F.
\end{align}
%The former equation represents the function $F=F(x,y)$ in an implicit form.
By putting $F= y^4 f$ and substituting it into the first 
equation of \eqref{quad2G30}, we obtain that the function $f$ 
is subject to the equation
\begin{align}\label{f10}
x^2 y^4 f^2 + 2x (y^2 - 2\ddd)  f + 1 = 0.
\end{align}
%Also, by putting $G = y^2 g$ and substituting it to 

From $F = y^4 f$, it is immediate to obtain 
\begin{align}\label{F2d}
y^4 \set F_{xx},\quad y^3 \set F_{xy}
\qandq y^2 \set F_{yy}.
\end{align}
These mean $y^6\set \Hess(F)$.
From the second equation of \eqref{quad2G30}, we obtain
\begin{align}\label{G11}
G = y^2(1 + xy^2 f),
\end{align}
and from this, we readily have
\begin{align}\label{G2d}
y^4\set G_{xx},\quad
y^3 \set G_{xy}\qandq
y \not|\, G_{yy}.
\end{align}
These mean $y^4\set \Hess(G)$.
Also, from \eqref{F2d} and \eqref{G2d}, we readily obtain 
$y^4\set K$.
In order to show $y^5 \not|\,K$,  differentiating \eqref{f10} by $x$,
we obtain 
$$
y^4 (xf^2 + x^2 ff_x) + (y^2 - 2 \ddd) (f + xf_x) = 0.
$$
As $\ddd\neq 0$, this means $y^4 \set (f + xf_x) $.
Writing $f = u(x) + y v(x,y)$ as before, this implies 
\begin{align}\label{fxx2}
f_{xx} = \frac{2c}{x^3} + yv_{xx}
\end{align}
for some constant $c$. If $c=0$, we obtain $y\set f$,
which cannot happen since no hyperplane section of $S$ can contain 
$5l$.
Hence from \eqref{fxx2}, we obtain $y\not|\, f_{xx}$.
So $y^5\not|\, y^4 f_{xx} = F_{xx}$.
As $y \not|\, G_{yy}$ as in \eqref{G2d}, we conclude $y^5 \not|\, F_{xx}G_{yy}$.
This means $y^5\not|\,K$.
Therefore, $y^4\set \bf H$ and $y^5 \not|\, \bf H$ hold.
Hence, {\em the divisor $\bm D_1$ includes 
the component $\pi_1\inv(l)$ with multiplicity precisely four.}
The discriminant of the quadratic polynomial ${\bf H}/y^4$ is given by
\begin{align}\label{Hess20}
\left(\frac{K}{y^4}\right)^2
-
4\frac{\Hess(F)}{y^4}\frac{\Hess(G)}{y^4}. 
\end{align}
Since  $y^5 \not|\,K$ as above,
$K/y^4$ is not divisible by $y$.
So $(K/y^4)^2$ is also not divisible by $y$.
On the other hand, we have $y^6\set \Hess(F)$ as above.
These mean that  the discriminant \eqref{Hess20} is not
divisible by $y$.
Therefore, {\em the line $l$ is not a branch divisor 
of the generically double covering $\bm D_S\to S$.}

\subsection{} Finally in this section, we consider a Segre surface $S$ 
whose symbol is [2(21)].
The normal forms of the defining equations of $S$ are given by 
$$
2\aaa X_0 X_1 +  X_1^2 + 2\bbb X_2 X_3 + X_3 ^2 + \bbb X_4^2=
2 X_0 X_1 + 2X_2 X_3 +  X_4^2 = 0,
$$
where $\aaa\neq \bbb$.
Singularities of $S$ are the point $e_0$, an $A_1$-point,
and $e_2$, an $A_3$-point.
The  line $l: = \{X_1 = X_3 = X_4 =0\}$ passes these singularities,
and it is contained in $S$.
The common tangent plane along $l\minus\{e_0,e_2\}$ is given by 
\begin{align}\label{ctp11}
T=\{X_1 = X_3 =0\}.
\end{align}

This time, we work on the affine open subset $\{X_2\neq 0\}$ in $\CP_4$
and use $x_i=X_i/X_2$ as coordinates on it.
The $A_2$-point $e_2$ is the origin.
From \eqref{ctp11}, we may use $(x,y):=(x_0,x_4)$ as coordinates
on $S$ around points of $l\minus\{e_0,e_2\}$. %, except the origin $e_2$.
On $S$, we again have $l = \{y=0\}$.
The second equation for $S$ becomes
\begin{align}\label{G12}
2 x_3 = - ( 2x  x_1 + y^2).
\end{align}
Substituting this into an affine form of the first equation of $S$,
we obtain 
\begin{align}\label{F12}
4(1+x^2) x_1^2 + 4 x (2\ddd + y^2)  x_1 + y^4  = 0.
\end{align}
If we let $x_1=0$ in this equation, we obtain $y^4=0$.
This means that the hyperplane section $\{X_1 = 0\}|_S$
is exactly $4l$.
The hyperplane $\{X_1=0\}$ is the unique one 
given in Lemma \ref{lemma:line4}.
We put
\begin{align}\label{FG}
F = x_1,\quad G = -2x_3,
\end{align}
and define $f$ and $g$ by $F = y^4 f$ and $G = y^2 g$.
From \eqref{F12}, we obtain that the function $f$ is subject to the equation
\begin{align}\label{quad55}
4(1+x^2) y^4 f^2 + 4x (2\ddd + y^2) f + 1 = 0.
\end{align}
Also, from \eqref{G12}, we have 
$$G = y^2( 2xy^2f + 1).$$
From this, we readily obtain 
\begin{align}\label{G2d2}
y^4 \set G_{xx},\quad 
y^3 \set G_{xy}\qandq
y\not|\, G_{yy}.
\end{align}
On the other hand, in the same way to the last case,
from $F=y^4f$, we obtain
\begin{align}\label{F2d2}
y^4 \set F_{xx},\quad y^3 \set F_{xy}
\qandq y^2 \set F_{yy}.
\end{align}
From \eqref{G2d2} and \eqref{F2d2}, we obtain
$$
y^6\set \Hess(F),\quad
y^4 \set \Hess(G)\qandq
y^4\set K.
$$
Also, differentiating \eqref{quad55} by $x$, we obtain
$$
2y^4\big\{xf^2 + (1+x^2) ff_x\big\} + (2\ddd+y^2) (f + xf_x) =0.
$$
Since $\ddd\neq 0$, this means $ y^4 \set (f + x f_x)$.
%From this we again obtain 
%\begin{align}\label{fxx3}
%f_{xx} = \frac{2c}{x^3} + yv_{xx}
%\end{align}
%for some constant $c\neq 0$.
%So in the same way to 
%the last case, we obtain $y^5 \not|\, F_{xx}G_{yy}$.
The rest is completely the same as the previous case
of $[(21)(11)]$, and we obtain that 
{\em the divisor $\bm D_1$ includes 
the component $\pi_1\inv(l)$ with multiplicity precisely four},
and that {\em the line $l$ is not a branch divisor 
of $\bm D_S\to S$.}

\medskip
Summing up the results in this section, as promised, we obtain the following conclusion.

\begin{proposition}\label{prop:D1c}
Let $S\subset\CP_4$ be a Segre surface and $l$ a line
which passes two singularities of $S$.
Then $l$ is not a branch divisor of the generically finite
double covering $\bm D_S\to S$.
\end{proposition}

Finally, we make a remark about multiplicities of the lines
as branch divisor of the generically finite
double covering $\bm D_S\to S$.
If a line on a Segre surface $S$ does not pass any singularity of $S$, then the multiplicity of the line is one
(Proposition \ref{prop:line1}).
If a line passes two singularities of $S$, it is not a branch
divisor (Proposition \ref{prop:D1c}).
When a line passes exactly one singularity of $S$,
we know from Proposition \ref{prop:line2} that 
the multiplicity of the line is at least two.
The exact multiplicities of these lines might be
determined by making the following heuristic argument solid.

Take a Segre surface $S$ whose symbol is $[1112]$ for example.
$S$ has exactly one singularity and 12 lines.
Precisely 4 of the lines pass the singularity. (See Table \ref{table1}.)
By realizing $S$ as a degeneration of smooth Segre surfaces,
each of the 4 lines can be obtained by unifying (or gathering) two lines
on smooth Segre surfaces through the degeneration.
By continuity, this would mean that the multiplicities of 
the 4 lines as branch divisor are precisely two.
Similarly, a Segre surface $S$ with symbol $[113]$ 
has exactly one singularity and  8 lines.
Precisely 4 of the lines pass the singularity.
Again by realizing $S$ as a degeneration of smooth Segre surfaces,
each of the 4 lines can be obtained by unifying (or gathering) three lines
on smooth Segre surfaces through the degeneration.
This would mean that the multiplicities of 
the 4 lines are precisely three.
In Table \ref{table1}, to each line which passes at least one singularity, the multiplicity of the line counted by using
degeneration from smooth Segre surfaces in this way is presented
in the $y$-column.
These multiplicities would be equal to the multiplicities as branch divisors
of the covering $\bm D_S\to S$.

\end{document}